\documentclass[11pt,a4paper,twoside]{article}
\usepackage{latexsym,amsfonts,amsmath, amsthm, amssymb}
\usepackage{fullpage}
\usepackage{xcolor,bm}
\usepackage[utf8x]{inputenc}
\usepackage{array}

\usepackage{graphicx}

\usepackage[normalem]{ulem}

\newtheorem{theorem}{Theorem}
\newtheorem{proposition}[theorem]{Proposition}
\newtheorem{lemma}[theorem]{Lemma}
\newtheorem{corollary}[theorem]{Corollary}

\theoremstyle{definition}

\newtheorem{definition}[theorem]{Definition}
\newtheorem{remark}[theorem]{Remark}

\newtheorem{note}[theorem]{Note}

\newcommand{\ls}{\leqslant}
\newcommand{\gs}{\geqslant}

\newcommand{\R}{\mathbb{R}}

\newcommand{\B}{\mathbf{B}}
\newcommand{\I}{\mathcal{I}}

\newcommand{\cZ}{\mathcal{Z}}
\newcommand{\PP}{\mathcal{P}}

\newcommand{\UU}{\mathbf{U}}

\newcommand{\II}{\mathfrak{I}}
\newcommand{\SU}{\mathcal{S}}
\newcommand{\U}{\mathcal{U}}

\newcommand{\norm}[1]{\left\lVert{#1}\right\rVert}
\newcommand{\abs}[1]{\left|{#1}\right|}

\DeclareMathOperator{\sign}{\mathrm{sign}}
\DeclareMathOperator{\intr}{\mathrm{int}}

\newcommand{\bfit}[1]{\textbf{\textit{#1}}}

\begin{document}

\title{\bf On the illumination of 1-symmetric convex bodies}

\medskip

\author{Wen Rui Sun and Beatrice-Helen Vritsiou}

\date{}

\maketitle

\begin{abstract}
\small
In \cite{Tikhomirov-2017} Tikhomirov verified the Hadwiger-Boltyanski Illumination Conjecture for the class of 1-symmetric convex bodies of sufficiently large dimension. We propose an alternative approach which allows us to settle the conjecture for this class in all dimensions in a uniform way. We also demonstrate that an alternative approach was indeed needed for the low dimensions. Finally, with this alternative method it is possible to solely use illuminating sets which consist of pairs of opposite directions; we thus also answer a question by Lassak, who has suggested this may be possible for any origin-symmetric convex body. As a consequence of this, we can also confirm the X-ray conjecture by Bezdek and Zamfirescu for all 1-symmetric convex bodies.
\end{abstract}

\section{Introduction}\label{sec:intro}

Let $K$ be a \emph{convex body} in the Euclidean space $\R^n$, by which we mean a convex, compact set with non-empty interior. Thinking of $K$ as a solid object, the notion of \emph{illuminating} it can be defined in a very intuitive way: given a boundary point $x$ of $K$, we say that a point `light source' $\bm{p}$ illuminates $x$ if the halfline originating from $\bm{p}$ and passing through $x$ intersects the interior of $K$ at a point not lying between $\bm{p}$ and $x$; in other words, if there exists $t>1$ such that $\bm{p}+t(x-\bm{p})\in \intr K$. By compactness, we can find finite sets $\{\bm{p}_1,\bm{p}_2,\ldots,\bm{p}_M\}$ of `light sources' such that each boundary point $x$ of $K$ is illuminated by at least one point in the set, and moreover by picking points $\bm{p}_i$ which are sufficiently distant from $K$, we can hope to find such sets of relatively small size $M$. The smallest cardinality that such a set of point `light sources' illuminating $K$ can have is called the \emph{illumination number} of $K$. We denote it by $\II(K)$.

A possibly less intuitive way to define illumination, which nevertheless leads to the same quantity for each body $K$, is the following: given a non-zero vector $d\in \R^n$ (a direction) and a boundary point $x$ of $K$, we say that $d$ illuminates $x$ if there exists $\varepsilon > 0$ such that $x+\varepsilon d\in \intr K$. A set of directions ${\cal D}= \{d_1, d_2,\ldots, d_{M^\prime}\}$ such that:
\begin{center}
for each boundary point $x$ of $K$, 
\\
there is at least one $d_i\in {\cal D}$ which illuminates $x$,
\end{center}
will be called an illuminating set for $K$. The smallest cardinality of such a set coincides with the illumination number $\II(K)$ of $K$, as defined previously. 

The former definition of illumination was introduced by Hadwiger \cite{Hadwiger-1960}, while the latter one was proposed by Boltyanski \cite{Boltyanski-1960}.

\bigskip

\noindent {\bf The Hadwiger-Boltyanski Illumination Conjecture.} For every convex body $K$ in $\R^n$, we should have $\II(K) \ls 2^n$.


Moreover, the inequality should be strict, except in the case of the cube and of its affine images (parallelepipeds) in $\R^n$.

\let\thefootnote\relax\footnote{\emph{Keywords:} illumination, symmetries of the cube, covering number, X-ray number, deep illumination}
\let\thefootnote\relax\footnote{\emph{2020 Mathematics Subject Classification:} 52A40, 52A37 (Primary); 52A20, 52C07 (Secondary)}

\medskip

An excellent reference on the history of the conjecture, related problems, and progress up to recent years is the survey \cite{Bezdek-Khan-survey}. We also refer to the monographs \cite{Bezdek-CMSmonograph, BMS-book} and the surveys \cite{Bezdek-1992, Boltyanski-Gohberg-1995, Martini-Soltan-1999}. It is worth noting that the illumination number $\II(K)$ of $K$ is always equal to its \emph{covering number}, that is, the minimum number of copies of $\intr K$ that we need in order to cover $K$. In other words, $\II(K) = N(K,\intr K)$, with the latter number having already appeared in the so-called Covering Conjecture. Levi in 1955 \cite{Levi-1955} had considered and fully settled the problem of bounding $N(K,\intr K)$ for planar convex bodies (showing that $N(K,\intr K) = 3$ for $K\subset \R^2$, except if $K$ is a parallelogram, in which case $N(K,\intr K) = 4$). Motivated by that, in 1957 Hadwiger \cite{Hadwiger-1957} posed the analogous question in higher dimensions.

\medskip

\noindent {\bf Hadwiger's Covering Problem.} Is it true that for every convex body $K$ in $\R^n$ one has that $N(K,\intr K)\ls 2^n$? Moreover, are parallelepipeds the only equality cases?

\medskip

There is also an equivalent formulation by Gohberg and Markus \cite{Gohberg-Markus-1960}, where one covers $K$ by smaller homothetic copies of it. In this paper we will exclusively work with Boltyanski's definition of illumination.

\bigskip

Aside from Levi's solution in $\R^2$, in all other dimensions the general problem is open. In dimension 3 Lassak \cite{Lassak-1984} has shown that, if $K$ is centrally symmetric, that is, $K=-K$ (or more generally, $K-x=x-K$ for some $x\in \R^3$), then $\II(K) \ls 8 = 2^3$. Thus, short of the equality cases, the conjecture in $\R^3$ is settled for symmetric convex bodies, but it remains open for the not-necessarily symmetric case, with the best known upper bound being $14$ (due to Prymak \cite{Prymak-2023}). We also refer to a very recent paper by Arman, Bondarenko and Prymak \cite{ABP-2024}, where the reader can find all the progress to date and the most recent improvements on the bounds for other low dimensions.

A longstanding general upper bound (which remains the best known when specialised to the symmetric case) was already given in 1964 by Erd\"{o}s and Rogers \cite{Erdos-Rogers-1964}:
\begin{equation*}
\II(K) = N(K,\intr K) \ls \frac{{\rm vol}(K-K)}{{\rm vol}(K)}\theta(K) \ls \frac{{\rm vol}(K-K)}{{\rm vol}(K)}n\bigl(\ln n + \ln\ln n + 5\bigr)
\end{equation*}
where $\theta(K)$ is the asymptotic lower density of the most economical covering of $\R^n$ by copies (translates) of $K$. Erd\"{o}s and Rogers adapted an earlier proof by Rogers \cite{Rogers-1957} which was giving the first polynomial-order, and essentially best known to date, bound on $\theta(K)$. Combining their conclusion with the celebrated Rogers-Shephard inequality \cite{Rogers-Shephard-1957}, one is led to the bound $\II(K) \ls C4^n \sqrt{n}\ln n$ for every convex body $K\subset \R^n$, where $C$ is an absolute constant (moreover, in the symmetric case one gets $\II(K) \ls C^\prime 2^n n\ln n$). More recently, subexponential improvements to this general upper bound were given in \cite{HSTV-2021}, \cite{CHMT-2022} and \cite{Galicer-Singer-2023}, with the latter two papers attaining almost exponential improvements. The main novelty in these three papers is the use of results from Asymptotic Convex Geometry on the concentration of volume in high-dimensional convex bodies. Note however that neither the initial approach in \cite{HSTV-2021}, nor the more recent refinements, can contribute anything to the bound in the symmetric case, which would be the most relevant one for the current paper.

\bigskip

The Illumination Conjecture has been fully settled for certain special classes of convex bodies. Again, we refer the reader to the survey \cite{Bezdek-Khan-survey} for a comprehensive list of references up to 2016. Just as examples, we mention that:
\begin{itemize}
\item Levi also showed in \cite{Levi-1955} that $\II(Q)=n+1$ for all smooth convex bodies $Q$ in $\R^n$ 
(\emph{smooth} here means that at each boundary point there is a unique supporting hyperplane).
\item Martini \cite{Martini-1985} settled the conjecture for the class of belt polytopes (which contains the zonotopes). This was later extended by Boltyanski and Soltan \cite{Boltyanski-Soltan-1990, Boltyasnki-Soltan-1992} to zonoids (limits of zonotopes), and by Boltyanski \cite{Boltyanski-1996} to belt bodies (see also \cite{Boltyanski-Martini-2001}; note that $n$-dimensional belt polytopes are already dense in the full space of $n$-dimensional convex bodies with respect to the Hausdorff metric, and that belt bodies are a `delicately-defined' proper extension of those, forming another dense subclass of convex bodies).

Observe that, just as the previous bullet point also shows, knowing the conjecture for a dense subclass of convex bodies is not enough to settle it for all bodies, and thus the aforementioned results used more intricate arguments to pass from zonotopes to zonoids and from belt polytopes to belt bodies.
\item The conjecture is fully settled for convex bodies of constant width. For dimensions $n\gs 16$, this is due to Schramm \cite{Schramm-1988}, who gave an exponential upper bound for the illumination number of such convex bodies which is slightly `corrected' by a small polynomial factor in the dimension; this `corrected' upper bound can be checked to be $< 2^n$ as soon as $n\gs 16$. For the remaining dimensions we have: \cite{Lassak-1997}, \cite{Weissbach-1996} (see also \cite[Section 11]{BLNP-2007}) dealing with $n=3$, \cite{Bezdek-Kiss-2009} dealing with $n=4$, and \cite{BPR-2022} dealing with $5\ls n\ls 15$.
\item Tikhomirov \cite{Tikhomirov-2017} settled the conjecture for 1-symmetric convex bodies of sufficiently large dimension (we recall the definition of 1-symmetric below). His result is the main motivation for the current paper.
\item Bezdek, Ivanov and Strachan \cite{BIS-2023} confirmed the conjecture for centrally symmetric cap bodies in dimensions $n=3$ (see also \cite{Ivanov-Strachan-2021}), $n=4$, and $n\gs 20$. They further showed that, if the cap body is 1-unconditional, then the Illumination Conjecture holds in all dimensions (and in that case $\II(K) \ls 4n$ once $n\gs 5$).
\item Gao, Martini, Wu and Zhang \cite{GMWZ-2024} verified the conjecture for polytopes which arise as the convex hull of the Minkowski sum of a finite subset of the lattice ${\mathbb Z}^n$ and of the unit-volume cube $\bigl[-\tfrac{1}{2},\,\tfrac{1}{2}\bigr]^n$.
\end{itemize}

Finally, Livshyts and Tikhomirov \cite{Livshyts-Tikhomirov-2020} settled the conjecture for convex bodies in sufficiently small neighbourhoods of the cube (with respect to either the geometric or the Hausdorff distance). Given that $\II(K) = N(K,\intr K)$ is an upper semicontinuous quantity (see e.g. \cite{Naszodi-2009}), the bound $2^n$ can already be deduced for bodies sufficiently close to $[-1,1]^n$, so their result is about settling the equality cases (and indeed they show that, if ${\rm dist}(K, [-1,1]^n)$ is small enough, and $K$ is not a parallelepiped, then $2^n -1$ is a \underline{sharp} upper bound for $\II(K)$).

\bigskip

Returning to the most relevant results for the current paper, let us first recall that a convex body $K$ in $\R^n$ is called \emph{$1$-symmetric} if
\begin{center}
$x=(x_1,x_2,\ldots,x_n)\in K\ $ implies that
$\ (\epsilon_1 x_{\sigma(1)}, \epsilon_2 x_{\sigma(2)},\ldots, \epsilon_n x_{\sigma(n)})\in K$ 
\\
for any choice of signs $\epsilon_i\in \{\pm 1\}$, $1\ls i\ls n$, and any permutation $\sigma$ on $[n]=\{1,2,\ldots, n\}$.
\end{center} 
$K$ is called \emph{$1$-unconditional} if we only require that 
\begin{center}
$x=(x_1,x_2,\ldots,x_n)\in K\ $ implies that
$\ (\epsilon_1 x_1, \epsilon_2 x_2,\ldots, \epsilon_n x_n\in K$ 
\\
for any choice of signs $\epsilon_i\in \{\pm 1\}$.
\end{center}
We already mentioned Tikhomirov's result, saying that there is an absolute constant $C_0$ such that, as long as the dimension $n\gs C_0$, then for every 1-symmetric convex body $K\subset \R^n$ which is not a multiple of the cube $[-1,1]^n$, one has $\II(K)\ls 2^n-1$ (no explicit lower bound for $C_0$ is given in \cite{Tikhomirov-2017}, nor have we tried to find one here; however, it should become clear in the next sections that a direct adaptation of Tikhomirov's method would not work in all low dimensions, and hence the restriction $n\gs C_0$ is not merely an artefact of the proof).

\smallskip

Of relevance, regarding the illumination of 1-symmetric or 1-unconditional convex bodies, are also the following three results in dimension 3: as we mentioned, Lassak \cite{Lassak-1984} showed that $\II(K)\ls 8$ for every origin-symmetric convex body $K$ in $\R^3$. Moreover, he showed this while using illuminating sets consisting of 4 pairs of opposite directions (and posed the question whether this is possible to do in higher dimensions as well). Bezdek \cite{Bezdek-1991} showed that $\II(P)\ls 8$ for any polytope in $\R^3$ which has a non-trivial affine symmetry. Finally Dekster \cite{Dekster-2000} obtained the same bound for any convex body $K$ in $\R^3$ which is symmetric about a plane. (Observe that any of these latter assumptions would be satisfied by a 1-symmetric or 1-unconditional body in $\R^3$.)

\medskip

The main purpose of this paper is to complement Tikhomirov's approach, and settle the illumination conjecture for 1-symmetric convex bodies in all dimensions (we deal with all dimensions $\gs 3$ in a uniform way, that is, our approach is not affected by the dimension, and thus we also provide an alternative proof for the cases which Tikhomirov had already settled). We also explain why a statement in Tikhomirov's paper (\cite[Lemma 8]{Tikhomirov-2017}) does not hold in all the generality in which it is stated (see Section \ref{sec:counterexamples}), and, simply for completeness, we give in its place a similar-looking result in the same spirit (see Theorem \ref{thm:illum-no-unit-squares}). That said, it should be stressed that Tikhomirov's result remains unaffected, as his proof can be `rectified in-house' by slightly rearranging the components of the proof, and where each one of them is used (see Remark \ref{rem:Tikh-proof-gap}). 

\smallskip

Finally, all the illuminating sets that we use in this paper can be chosen to consist of pairs of opposite directions. We thus answer Lassak's conjecture too for 1-symmetric convex bodies of dimension $n\gs 3$, and, as a consequence, we can also confirm the X-ray conjecture by Bezdek and Zamfirescu in this class (see Section \ref{sec:main-results} for details).

\medskip

To describe our results and approach in detail, and compare it with Tikhomirov's, we first need some additional notation and definitions, so we defer doing so to Subsection \ref{subsec:Tikh-and-our-method}. 

%

\section{Preliminaries, and precise statements of our results}\label{sec:prelims}

We write $[n]$ for the set $\{1,2,\ldots,n\}$, and $e_1, e_2,\ldots, e_n$ for the standard basis vectors of $\R^n$. For any vector $x\in \R^n$, we will denote by $\cZ_x$ the set $\{i\in [n]: x_i=0\}$. 

Given a subset $A$ of $\R^n$, we will denote its interior and its boundary by $\intr A$ and by ${\rm bd} A$ or $\partial A$ respectively. Recall that if $A$ is a non-empty convex set, then its affine hull
\begin{equation*}
{\rm aff}(A) := \bigl\{\mu_1a_1+\mu_2a_2+\cdots+\mu_\ell a_\ell\,:\,  \ell \gs 1,\, a_i\in A \ \hbox{and}\ \mu_i\in \R \ \hbox{with}\ \mu_1 + \mu_2 +\cdots + \mu_\ell = 1\bigr\}
\end{equation*}
coincides with the smallest affine subspace of $\R^n$ which contains $A$, and that, in the subspace topology on ${\rm aff}(A)$, $A$ has non-empty interior. We call this the relative interior of $A$ and denote it by ${\rm relint} A$. Moreover, we call $A\setminus {\rm relint} A$ the relative boundary of $A$, and denote it by ${\rm relbd}A$.

\medskip

Recall that $K$ is called a convex body if it is a compact convex set with non-empty interior. If $K$ is also origin-symmetric, that is, $K=-K$, then we know that $K$ is the unit ball of a certain norm on $\R^n$, which is given by $x\in \R^n \, \mapsto\, \|x\|_K:= \inf\{t>0: x\in tK\}$.

\smallskip

We also recall that the illumination number of any convex body is an affine invariant: namely $\II(K) = \II(TK+z)$ for any invertible linear transformation $T\in {\rm GL}(n)$ and any (translation) vector $z$. 

Therefore, without loss of generality, we can assume that all the 1-symmetric convex bodies $\B\subset \R^n$ which we consider satisfy $e_i\in \partial \B$ for all $i\in [n]$. This is equivalent to saying that $\|e_i\|_\B = 1$ for all $i\in [n]$. We will denote this subclass of $n$-dimensional 1-symmetric convex bodies by $\SU^n$ \emph{(observe that any $n$-dimensional $1$-symmetric convex body has a dilation which is a member of this subclass)}.

\smallskip

Analogously, by simply multiplying by a diagonal matrix (but not necessarily a multiple of the identity now), we can assume that a given 1-unconditional convex body $\UU\subset \R^n$ satisfies $\|e_i\|_\UU=1$ for all $i\in [n]$. We will denote this subclass of $n$-dimensional 1-unconditional convex bodies by $\U^n$.

\smallskip

We recall the following necessary and sufficient condition for illumination (even though in the sequel we will mostly work with the definition by Boltyanski mentioned in the Introduction).

\smallskip

\noindent {\bf Remark/Criterion A.}
Let $K$ be a convex body in $\R^n$, let $x\in \partial K$, and let $d$ be a direction in $\R^n$. Then $d$ illuminates $x$ if and only if,
\begin{center}
for every outer normal vector $u$ of $K$ at $x$, we have $\langle d, u\rangle < 0$
\end{center}
(where $\langle \cdot, \cdot \rangle$ denotes the standard dot product on $\R^n$).

\medskip

\begin{lemma}\label{lem:smaller-coordinates}
Let $\UU$ be a 1-unconditional convex body in $\R^n$. Suppose that $x$ is a point in $\UU$, and that $y\in \R^n$ satisfies:
\begin{equation*}
\hbox{for all $i\in [n]$},\ \ |y_i|\ls |x_i|.
\end{equation*}
Then $y\in \UU$ as well.

Moreover, if we have that
\begin{equation*}
\hbox{for all $i\in [n]$},\ \ |y_i| < |x_i| \ \hbox{or}\ |y_i|=|x_i|=0,
\end{equation*}
then $y\in \intr\UU$.
\end{lemma}
\begin{proof}
This follows quickly from the 1-unconditionality of $\UU$ and convexity: given these, $\UU$ contains the (possibly lower-dimensional) parallelepiped $\PP_x$ with vertices the point $x$ and all its coordinate reflections. It is clear that, if $|y_i|\ls |x_i|$ for all $i\in [n]$, then $y\in \PP_x\subseteq \UU$.

\smallskip

Moreover, if we actually had that $x\in \intr\UU$, then $\PP_x\subset \intr\UU$ (as the convex hull of interior points of $\UU$). Thus, if we now assume that, for all $i\in [n]$, $|y_i| < |x_i|$ or $|y_i| = |x_i| = 0$, then we can find some $t \in (0,1)$ such that $|y_i| \ls t |x_i|$ for all $i\in [n]$, and we can apply the above with the point $t x\in \intr\UU$ instead of $x$: we will have $y\in \PP_{t x}\subset \intr\UU$, as claimed.
\end{proof}

\begin{corollary}\label{cor:uncond-illum}
Let $\UU$ be a 1-unconditional convex body in $\R^n$, and let $x\in \partial \UU$. Then $x$ is illuminated by any direction $d\in \R^n$ which satisfies
\begin{center}
$\cZ_d = \cZ_x$, and $d_i\cdot x_i < 0$ for all $i\in [n]\setminus \cZ_x$
\end{center}
(recall that $\cZ_x$ is the set of indices in $[n]$ which correspond to the zero coordinates of $x$). 

\smallskip

In particular, $\UU$ is illuminated by the set $\{-1,0,1\}^n\setminus \{\vec{0}\}$. 
\end{corollary}
\begin{proof}
Simply note that there is $\varepsilon >0$ small enough so that, for the displaced vector $x+\varepsilon d$, we will have $|(x+\varepsilon d)_i| < |x_i|$ for all $i\in [n]\setminus\cZ_x$. Thus $x+\varepsilon d\in \intr\UU$ by the previous lemma.
\end{proof}

It is also worth clearly stating the following

\begin{remark}\label{rem:sym-and-uncond-in-cube}
If $\B$ is a 1-symmetric or 1-unconditional convex body in $\R^n$, and $x\in \B$, then, by Lemma \ref{lem:smaller-coordinates}, we also have that $|x_i|e_i\in \B$ for all $i\in [n]$. Thus, if $\B\in \SU^n$ or $\in \U^n$, then $\|x\|_\infty:=\max\limits_{i\in [n]}|x_i|\ls 1$. In other words, $\B\subseteq [-1,1]^n$.
 \end{remark}

\smallskip

\begin{lemma}\label{lem:unique-max-coord}
Let $\B$ be a 1-symmetric convex body in $\R^n$, and let $x\in \partial\B$. Suppose that there is a unique index $i_0\in [n]$ such that $|x_{i_0}| = \norm{x}_\infty$. Then $x$ can be illuminated by $-\sign(x_{i_0})e_{i_0}$.
\end{lemma}
\begin{note}
This is used in \cite{Tikhomirov-2017} as well (see the first part of the proof of \cite[Lemma 8]{Tikhomirov-2017}). There it is proved using Criterion A for illumination; here we provide one more proof, this one based on the classical definition by Boltyanski.
\end{note}
\begin{proof}
Let us look at the coordinates of the displaced vector $\widetilde{x}=x+ |x_{i_0}|(-\sign(x_{i_0})e_{i_0})$: its $i_0$-th coordinate is zero, while, for $i\neq i_0$, $\widetilde{x}_i=x_i$. 

Since $\B$ is 1-symmetric, the following convex combination is contained in it:
\begin{align*}
w_x:=\,&\frac{1}{n}\left[\sum_{i=1}^n|x_i|e_i\  +\ \sum_{i\neq i_0}\left(|x_i|e_{i_0} + |x_{i_0}|e_i \, +\, \sum_{j\notin\{i,i_0\}}|x_j|e_j\right)\right]
\\
=\,&\left(\frac{1}{n}\sum_{j=1}^n|x_j|\right)e_{i_0} \ +\  \sum_{i\neq i_0}\left(\frac{n-1}{n}|x_i| + \frac{1}{n}|x_{i_0}|\right)e_i.
\end{align*}
We can observe that
\begin{equation*}
\frac{1}{n}\sum_{j=1}^n|x_j| \gs \frac{1}{n}|x_{i_0}| > 0,\ \ \hbox{and also}\ \frac{n-1}{n}|x_i| + \frac{1}{n}|x_{i_0}| > \frac{n-1}{n}|x_i| + \frac{1}{n}|x_i| = |x_i| \ \hbox{for all}\ i\neq i_0,
\end{equation*}
since $|x_{i_0}|$ is the unique maximum coordinate of $x$ in absolute value. Thus we can apply Lemma \ref{lem:smaller-coordinates} with $w_x$ and $\widetilde{x}$ to conclude that $\widetilde{x}\in \intr\B$.
\end{proof}

\medskip

An equivalent statement for the next lemma, which will be given right after, can be viewed in some sense as a more broadly applicable version of the second part of Lemma \ref{lem:smaller-coordinates}.

\begin{lemma}\label{lem:affine-set}
Let $K$ be a convex body in $\R^n$, and let $H$ be an affine subspace of $\R^n$. Suppose that $(\intr K) \cap H \neq \emptyset$. Then 
\begin{center}
${\rm relint}(K\cap H) = (\intr K) \cap H\,$ and $\,{\rm relbd}(K\cap H) = (\partial K)\cap H$.
\end{center}
\end{lemma}
\begin{proof}
From our assumptions, $K\cap H$ is a non-empty convex set. It is clear that $(\intr K) \cap H\subseteq {\rm relint}(K\cap H)$ and ${\rm relbd}(K\cap H) \subseteq (\partial K)\cap H$. Observe also that $H= {\rm aff}(K\cap H)$ because $(\intr K) \cap H \neq \emptyset$.

Fix $x_0\in (\intr K) \cap H$. Also, let $y\in {\rm relint}(K\cap H)$. Note that, for every $s\in \R$, we will have that
\begin{equation*}
y+s(y-x_0)\in {\rm aff}(H) = H = {\rm aff}(K\cap H).
\end{equation*}
Since $y\in {\rm relint}(K\cap H)$, we can find $s_y>0$ small enough so that $y+s_y(y-x_0)\in K\cap H$. Then we can write
\begin{equation*}
y = \frac{s_y}{1+s_y}x_0 \,+\,\frac{1}{1+s_y}\bigl(y+s_y(y-x_0)\bigr),
\end{equation*}
which shows that $y\in \intr K$, as a non-trivial convex combination of two points in $K$ with at least one of them being an interior point. Thus $y\in (\intr K) \cap H$, which implies that ${\rm relint}(K\cap H)\subseteq (\intr K) \cap H$ as well. 

\smallskip

From this it also follows that $(\partial K)\cap H\subseteq {\rm relbd}(K\cap H)$, which completes the proof.
\end{proof}

As mentioned, we are going to use an equivalent statement of Lemma \ref{lem:affine-set}, more obviously relevant to illumination. In some sense, this could be viewed as a counterpart to Lemma \ref{lem:smaller-coordinates} and Corollary \ref{cor:uncond-illum}, that now applies to arbitrary convex bodies.

Similarly to before, note that, if $p\in {\rm relbd}(K\cap H)$, and if $d^\prime$ is a non-zero vector in $H-p = \{u-p:u\in H\}$, then $p+sd^\prime\in H$ for any $s\in \R$ (this is because $H-p$ is a linear subspace of $\R^n$, and thus $d^\prime\in H-p \,\Rightarrow\,sd^\prime\in H-p$).

\smallskip

\begin{corollary}\label{cor:affine-set}
Given the same assumptions as before, if $p\in {\rm relbd}(K\cap H)$, and $d^\prime$ is a non-zero vector in $H-p\subset \R^n$ such that $p+\varepsilon d^\prime \in {\rm relint}(K\cap H)$ for some $\varepsilon > 0$, then $p+\varepsilon d^\prime\in \intr K$.

In other words, if $p$ is $(K\cap H)$-illuminated by $d^\prime$, within $H= {\rm aff}(K\cap H)$, then it is also $K$-illuminated by $d^\prime$, viewed within $\R^n$ now.
\end{corollary}
\begin{proof}
By Lemma \ref{lem:affine-set}, we know that ${\rm relint}(K\cap H) = (\intr K) \cap H$, and hence
\begin{equation*}
p+\varepsilon d^\prime \in {\rm relint}(K\cap H) = (\intr K) \cap H\ \ \Rightarrow\ \ p+\varepsilon d^\prime\in \intr K,
\end{equation*}
as claimed.
\end{proof}

As we will see in the sequel, to justify the illumination of certain boundary points, it makes sense to also resort to the somewhat indirect reasoning allowed by Corollary \ref{cor:affine-set}: for instance, in certain situations this will enable us to employ partially inductive arguments, or facilitate working with directions which have some zero coordinates.

\medskip

We finish the preliminaries with a basic remark on the illumination of certain 1-unconditional convex bodies. For notational convenience, we write ${\bm 1} = e_1+e_2+\cdots +e_n$.

\begin{lemma}\label{lem:small-distance-to-cube}
Let $\UU$ be a 1-unconditional convex body in $\U^n$ (recall that this means that $\|e_i\|_\UU=1$ for all $i\in [n]$). Assume that $\frac{1}{2}{\bm 1}\in \intr\UU$. 
Then $\U$ is illuminated by the set $\{-1,1\}^n$.
\end{lemma}
\begin{proof}
Let $x\in \partial \UU$, and pick a direction $d_x$ from $\{-1,1\}^n$ such that $d_{x,i}\cdot x_i < 0$ whenever $x_i\neq 0$. Then the displaced vector $x+ \frac{1}{2}d_x$ satisfies
\begin{equation*}
\big|(x+ \tfrac{1}{2}d_x)_i\big| = \max\bigl(|x_i|-\tfrac{1}{2},\,\tfrac{1}{2}-|x_i|\bigr)\ls \max\bigl(1-\tfrac{1}{2},\,\tfrac{1}{2}-|x_i|\bigr)=\tfrac{1}{2}
\end{equation*}
whenever $x_i\neq 0$, and $\big|(x+ \tfrac{1}{2}d_x)_i\big| = \frac{1}{2}$ if $x_i=0$ (recall that $\UU\in \U^n$ implies that $\UU\subset [-1,1]^n$, hence $|x_i|\ls 1$ for all $i\in [n]$). In short, all coordinates of $x+ \frac{1}{2}d_x$ are $\ls \frac{1}{2}$ in absolute value. 

We now observe that, by our main assumption, there is $\eta>0$ such that
\begin{equation*}
\bigl(\tfrac{1}{2}+\eta\bigr){\bm 1} \in \UU.
\end{equation*}
We can thus apply Lemma \ref{lem:smaller-coordinates} to conclude that $x+ \frac{1}{2}d_x\in \intr\UU$, and that $d_x$ illuminates $x$.
\end{proof}

Note that this result is related to the upper semicontinuity of the illumination number $\II(\cdot)$ at the `point' $[-1,1]^n$; see \cite[Proposition 2.2]{Naszodi-2009} for a general argument using the language of covering.

\subsection{Comparing with Tikhomirov's approach}\label{subsec:Tikh-and-our-method}

Tikhomirov's approach has as its starting point the conclusion of Corollary \ref{cor:uncond-illum}: since $\{-1,0,1\}^n\setminus\{\vec{0}\}$ is an illuminating set for any 1-symmetric convex body $\B$ in $\R^n$, the objective is to look among its different subsets for `efficient' illuminating sets of 1-symmetric convex bodies, and to verify that, for each $\B\in \SU^n$, there is at least one such subset that works.

\smallskip

Tikhomirov distinguishes two main cases, and treats the 1-symmetric convex bodies in each case differently, based on their geometric distance from the cube $[-1,1]^n$. Given that, if $\B\in \SU^n$, then $\B\subseteq [-1,1]^n$, and the vectors $e_i$ are common points of $\B$ and of $[-1,1]^n$, we can simply define the geometric distance here as follows:
\begin{equation*}
{\rm dist}(\B,[-1,1]^n) := \inf\bigl\{\alpha \gs 1: [-1,1]^n\subseteq \alpha\B\bigr\}.
\end{equation*}
We can also check that, in this setting, ${\rm dist}(\B,[-1,1]^n) = \|e_1+e_2+\cdots+e_n\|_\B^{-1} = \|{\bm 1}\|_\B^{-1}$.

\medskip

Tikhomirov's result follows from the following two propositions.

\smallskip

\noindent {\bf Proposition B.} \emph{(\cite[Proposition 5]{Tikhomirov-2017})} Let $n\gs 2$ and $\B\in \SU^n$ with $1< {\rm dist}(\B,[-1,1]^n) < 2$. Then one of the following subsets of $\{-1,0,1\}^n\setminus\{\vec{0}\}$ is an illuminating set for $\B$:
\begin{itemize}
\item $T_1:= \bigl\{(\epsilon_1,\epsilon_2,\ldots,\epsilon_n)\in \{-1,1\}^n:\exists\,i\ls n-1\ \hbox{with}\ \epsilon_i=-1\bigr\} \cup \{e_1+e_2+\cdots +e_{n-1}\}$.
\item $T_2:= \bigl(\{-1,1\}^{n-1}\times\{0\}\bigr)\cup\{\pm e_n\}$.
\end{itemize}

\smallskip

\noindent {\bf Proposition C.} \emph{(essentially \cite[Proposition 10]{Tikhomirov-2017})} Let $n\gs 2$ and $\B\in \SU^n$ with ${\rm dist}(\B,[-1,1]^n) \gs 2$. Then $\B$ can be illuminated by a set $T_3$ of the form:
\begin{equation*}
T_3 = \bigl(\{-1,1\}^{n-1}\times\{0\}\bigr) \cup R_0
\end{equation*}
where $R_0$ is any subset of $\{-1,0,1\}^n\setminus\{\vec{0}\}$ with the property:
\begin{gather*}
\hbox{for every $k\ls \lceil\frac{n}{2}\rceil$, and for every $y\in \{-1,0,1\}^n$ with exactly $k$ non-zero coordinates,}
\\ 
\hbox{there exists $z=z_y\in R_0$ with exactly $2k-1$ non-zero coordinates} \tag{P1}
\\ 
\hbox{and such that $y_i=z_i$ for all $i\in [n]$ for which $y_i\neq 0$.}
\end{gather*}

\smallskip

The remaining key ingredient in Tikhomirov's method is a probabilistic argument that allows him to prove that, once the dimension $n$ gets sufficiently large, then there exist subsets $R_0$ of $\{-1,0,1\}^n\setminus\{\vec{0}\}$ with property (P1) which also have cardinality $\abs{R_0}\ls \frac{2^n}{n}$. Thus, in sufficiently high dimensions, we can pick sets $T_3$ of the form described above which are `efficient' illuminating sets for any $\B\in \SU^n$ with ${\rm dist}(\B,[-1,1]^n) \gs 2$.

\bigskip

Observe that we could have also distinguished two main cases in a slightly different way. For $\B\in \SU^n$ define
\begin{equation*}
m_\B:=\max\bigl\{k\in \{1,2,\ldots,n\}: e_1+e_2+\cdots + e_k\in \B\bigr\}. 
\end{equation*}
Equivalently, $m_\B$ is the largest dimension of a unit subcube $[-1,1]^k\times\{0\}^{n-k}$ contained in $\B$. By our chosen normalisation for $\B$ (that is, given that $\B\in \SU^n$), we have that $m_\B$ is certainly $\gs 1$ (and thus, well-defined), and also that ${\rm dist}(\B,[-1,1]^n) = 1$ if and only if $m_\B=n$.

\smallskip

Note now that ${\rm dist}(\B,[-1,1]^n) \gs 2$ implies that $m_\B\ls \frac{n}{2}$. Indeed, assume that we have $m_\B > \frac{n}{2}$, and also initially suppose that $n$ is odd. Then by the 1-symmetry of $\B$, and Lemma \ref{lem:smaller-coordinates}, we will have that
\begin{center}
both $e_1+e_2+\cdots +e_{\frac{n+1}{2}}\in \B$, and $e_1+e_{\frac{n+3}{2}}+e_{\frac{n+5}{2}}+\cdots + e_n\in \B$.
\end{center}
Thus we will also have that $e_1+\frac{1}{2}\bigl(e_2+e_3+\cdots +e_n\bigr)\in \B$. By the 1-symmetry again, we get that
\begin{equation*}
\frac{n+1}{2n} {\bm 1 } = \frac{1}{n}\sum_{i\in [n]}\left(e_i+\frac{1}{2}\sum_{j\neq i}e_j\right) \in \B,
\end{equation*}
which shows that ${\rm dist}(\B,[-1,1]^n) = \|{\bm 1}\|_\B^{-1} \ls \frac{2n}{n+1} < 2$.

The proof is analogous and easier in the case that $n$ is even.

\begin{remark}\label{rem:Tikh-proof-without-dist}
The proof of \cite[Proposition 10]{Tikhomirov-2017}/Proposition C can be adjusted just slightly so that the assumption that ${\rm dist}(\B,[-1,1]^n) \gs 2$ is only used in order to ensure that $m_\B\ls \frac{n}{2}$. From then on, one can reach the same conclusion without using the distance assumption anymore. 

Indeed, if we know that $m_\B\ls \frac{n}{2}$, and we first consider, as Tikhomirov does, a boundary point $x$ such that
\begin{equation}\label{eqp:rem:Tikh-proof-without-dist}
\abs{\{i\in [n]: |x_i|=\|x\|_\infty\}} > \lceil\tfrac{n}{2}\rceil,
\end{equation}
then we can observe that $\|x\|_\infty$ will have to be $\in (0,1)$. This is given the normalisation $\B\in \SU^n$, which further implies that $\pm e_n$ cannot be an outer normal vector of $\B$ at $x$. But then, as Tikhomirov shows, one of the directions in $\{-1,1\}^{n-1}\times \{0\}$ will illuminate $x$, because of \eqref{eqp:rem:Tikh-proof-without-dist} as well.

\smallskip

For any of the remaining boundary points $y$, we will have that $\abs{\{i\in [n]: |y_i|=\|y\|_\infty\}} \ls \lceil\tfrac{n}{2}\rceil$, and then, exactly as Tikhomirov shows, a direction from the set $R_0$ with property (P1) will illuminate $y$ (no matter which such set $R_0$ we chose to work with).
\end{remark}

\begin{remark}\label{rem:Tikh-proof-gap}
Tikhomirov breaks down the proof of \cite[Proposition 5]{Tikhomirov-2017}/Proposition B into the proof of two lemmas. 

\smallskip

\cite[Lemma 7]{Tikhomirov-2017} states that, if $\B\in \SU^n$ and $1< {\rm dist}(\B,[-1,1]^n) < 2$, then either $\B$ is illuminated by the set $T_1$ (using the notation of Proposition B), or we necessarily have that $\|e_1+e_2\|_\B > 1$ (equivalently, we have $e_1+e_2\notin\B$). Note that of course both statements can be true as well.

\medskip

On the other hand, \cite[Lemma 8]{Tikhomirov-2017} states that, if $n\gs 2$ and $\B\in \SU^n$ satisfies $\|e_1+e_2\|_\B > 1$, then $\B$ can be illuminated by the set $T_2$ (again, using the notation of Proposition B). In the next section we will show that {\bf this statement does NOT hold in such generality}.

\medskip

We should stress however that \cite[Lemma 8]{Tikhomirov-2017} could be completely bypassed in Tikhomirov's approach. Indeed, if we have that $\|e_1+e_2\|_\B > 1\,\Leftrightarrow\,e_1+e_2\notin \B$, then we will certainly have that $m_\B\ls \frac{n}{2}$. Thus, Tikhomirov's method works based on the following main ingredients:
\begin{itemize}
\item If $\B\in \SU^n$ satisfies $1< {\rm dist}(\B,[-1,1]^n) < 2$ and $\|e_1+e_2\|_\B = 1$, then necessarily $\B$ is illuminated by the set $T_1$ and $\II(B)\ls 2^n -1$.
\item On the other hand, if ${\rm dist}(\B,[-1,1]^n) \gs 2$ or $\|e_1+e_2\|_\B > 1$, then we certainly have that $m_\B\ls \frac{n}{2}$. Thus, as we explained in Remark \ref{rem:Tikh-proof-without-dist}, a slight adjustment of the proof of \cite[Proposition 10]{Tikhomirov-2017} gives that $\B$ can be illuminated by some set $T_3$ (of the form described in Proposition C).
\item The final part of Tikhomirov's proof, which is probabilistic, guarantees the existence of sets $T_3$ with cardinality $\abs{T_3}\ls 2^{n-1} + \frac{2^n}{n}$, as long as $n$ is sufficiently large.
\end{itemize}
\end{remark}

\medskip

Our approach has some parallels to Tikhomirov's approach. We adopt a similar starting point, but instead of working with illuminating sets which are subsets of $\{-1,0,1\}^n\setminus\{\vec{0}\}$, we will work with subsets of the set
\begin{equation}\label{def:big-illum-set}
    G^n(\delta) := \left\{d\in \R^n: \exists\, i\in [n]\ \hbox{such that}\ d=\pm e_i\,+\!\sum_{j \in [n] \backslash \{i\}} 
   \! \pm \delta e_j \right\},
\end{equation}
(where $\delta\in (0,1)$ will be a small, suitably chosen parameter), or very minor variations of such subsets. Note that the main feature of all directions in $G^n(\delta)$ is that they are small perturbations of some standard basis vector $e_i$ (`small' because, as we will see, $\delta$ will be chosen smaller and smaller as the dimension grows).

\medskip

Our main results are as follows.

\smallskip

\noindent {\bf Theorem D.} For every $n\gs 2$, we can find a subset $\I^n$ of $G^n\bigl((n+1)^{-1}\bigr)$ with cardinality $\abs{\I^n}=2^n$ such that $\I^n$ will illuminate all 1-symmetric convex bodies in $\R^n$.

\medskip

\noindent {\bf Theorem E.} \emph{(see Theorem \ref{thm:illum-everything-but-the-cube})} Let $n\gs 3$ and let $\B\in \SU^n$ with ${\rm dist}(\B,[-1,1]^n) > 1$. Then we can slightly modify the set $\I^n$ from Theorem D to get an illuminating set for $\B$ as follows: we completely remove one pair of directions from $\I^n$, and also slightly alter the `small' coordinates of another pair (the altered directions will still be perturbations of standard basis vectors; also, the directions that we choose to replace/modify will have some relation to the pair of removed directions so as to `make up' for the absence of those too). It follows that $\II(\B) \ls 2^n -2$.

\medskip

Theorem E settles the Illumination Conjecture for all 1-symmetric convex bodies of any dimension $n\gs 3$. Moreover, as we will see, we can essentially use the same illuminating set for all 1-symmetric convex bodies in $\R^n$ which are not multiples of the cube, with the only variation coming from how small the `small' coordinates of the pair of altered directions will be (and this will only depend on how small $\,{\rm dist}(\B,[-1,1]^n) - 1\,$ is).

\medskip

Just for completeness, we can also provide some counterparts for the claimed \cite[Lemma 8]{Tikhomirov-2017} and for \cite[Proposition 10]{Tikhomirov-2017}, the latter dealing with 1-symmetric convex bodies $\B$ which satisfy ${\rm dist}(\B,[-1,1]^n) \gs 2$; see Theorems \ref{thm:illum-no-unit-squares} and \ref{thm:illum-small-unit-subcubes}. We call these results counterparts because they deal with the same subclasses of 1-symmetric convex bodies as \cite[Lemma 8 and Proposition 10]{Tikhomirov-2017} respectively, and because we can also use illuminating sets which are similar in form and function to $T_2$ and $T_3$.

The rest of the paper is organised as follows. In Section \ref{sec:counterexamples} we elaborate on why \cite[Lemma 8]{Tikhomirov-2017} cannot hold in such generality, and also on why Tikhomirov's core approach may not be suitable for `efficient' illumination of all cases in low dimensions (even if the proofs were to be modified). We prove Theorem D in Section \ref{sec:main-illuminating-set}. Then Theorem E/Theorem \ref{thm:illum-everything-but-the-cube}, as well as Theorems \ref{thm:illum-no-unit-squares} and \ref{thm:illum-small-unit-subcubes}, are established in Section \ref{sec:main-results}. Finally, in Section \ref{sec:almost-1-sym} we extend our main result to a slightly larger subclass, of `almost' 1-symmetric convex bodies, or $(1+\delta_n)$-symmetric bodies, where $\delta_n$ will depend on the dimension (for the exact definition, see Section \ref{sec:almost-1-sym}).

\section{`Tricky' convex bodies}\label{sec:counterexamples}

We present examples of 1-symmetric convex bodies in $\mathcal{S}^n$ which cannot be illuminated by the corresponding set $T_2$ in $\R^n$ (recall the notation in Proposition B) or even small variations of it (even though they satisfy the assumptions of \cite[Lemma 8]{Tikhomirov-2017}, which we recalled in Remark \ref{rem:Tikh-proof-gap}).

\smallskip

\begin{enumerate}
\item[1a.] Consider the body $\B_1^3\subset \R^3$ whose vertices are all the coordinate reflections and permutations of the points $e_1$ and $(\frac{1}{2},\frac{1}{2},\frac{1}{2})$. We can show that not only doesn't $T_2$ illuminate $\B_1^3$, but that in general any minor tweaks of the method in \cite{Tikhomirov-2017} will fail to work for $\B_1^3$: indeed, we will see that no subset of $\{-1,0,1\}^3\setminus\{\vec{0}\}$ with fewer than 10 elements can illuminate $\B_1^3$.

Note that among the outer normals at $e_1$ we have the vector $(1,0,0)$ and the four vectors $(1,\pm 1,0), (1,0,\pm 1)$. Because of the outer normal $(1,0,0)$, any direction $d$ which will illuminate $e_1\in \partial \B_1^3$ must satisfy $d_1 < 0$ (recall Criterion A). Thus, if $d$ is taken from $\{-1,0,1\}^3$, we must have $d_1 = -1$. On the other hand, because of the latter group of outer normals, any direction $d\in \{-1,0,1\}^3\setminus\{\vec{0}\}$ which will illuminate $e_1\in \partial \B_1^3$ cannot have non-zero 2nd or 3rd coordinates, because if it did, then $d$ would be orthogonal to at least one of the four vectors $(1,\pm 1,0), (1,0,\pm 1)$. In conclusion only the element $-e_1$ of $\{-1,0,1\}^3\setminus\{\vec{0}\}$ is an illuminating direction for $e_1\in \partial \B_1^3$. 

Similarly we see that we must include all the directions of the form $\pm e_i, i\in [3]$, to illuminate all the coordinate permutations and coordinate reflections of $e_1$, which are also vertices of $\B_1^3$. 

On the other hand, among the outer normals at $(\frac{1}{2},\frac{1}{2},\frac{1}{2})$ are the vectors $(1,1,0), (1,0,1)$ and $(0,1,1)$. Thus, none of the directions $\pm e_i, i\in [3]$, can illuminate $(\frac{1}{2},\frac{1}{2},\frac{1}{2})$, and we need instead a direction from $\{-1,0,1\}^3$ which has at least two non-zero (in fact, negative) coordinates. If we pick for this purpose the direction $(-1,-1,-1)$, then it will illuminate $(\frac{1}{2},\frac{1}{2},\frac{1}{2})$, but it won't illuminate any of its (non-trivial) coordinate reflections. If we instead pick a direction such as $(-1,-1,0)$, then we would be able to illuminate both $(\frac{1}{2},\frac{1}{2},\frac{1}{2})$ and $(\frac{1}{2},\frac{1}{2},\,-\!\!\!\!\!-\frac{1}{2})$. 

In this way, we can check that we need at least 4 directions from $\{-1,0,1\}^3\setminus\{e_i: i\in [3]\}$ to illuminate $(\frac{1}{2},\frac{1}{2},\frac{1}{2})$ and all its coordinate reflections.
\item[1b.] Consider the set $\B_1^4\subset \R^4$ whose vertices are all the coordinate reflections and permutations of the points $e_1$ and $(\frac{1}{2},\frac{1}{2},\frac{1}{2},\frac{1}{2})$. An analogous argument shows that $\B_1^4$ cannot be illuminated by any subset of $\{-1,0,1\}^4\setminus\{\vec{0}\}$ which contains fewer than 16 elements (and thus it certainly cannot be illuminated by $T_2$).
\item[1c.] In general, consider the set $\B_1^n\subset \R^n$ whose vertices are all the coordinate reflections and permutations of the points $e_1$ and $\frac{1}{2}(e_1+e_2+\cdots+e_n)$. We can similarly check that $\B_1^n$ cannot be illuminated by the set $T_2\in \R^n$, and at the very least we would need to enlarge $T_2$ to the set
\begin{equation*}
T_2^\prime := \bigl(\{-1,1\}^{n-1}\times\{0\}\bigr)\cup\{\pm e_i: i\in [n]\}
\end{equation*}
which contains $2^{n-1} + 2n$ directions.
\end{enumerate}

The last point might perhaps suggest that we could replace the set $T_2$ by $T_2^\prime$, and that this potentially is enough  in dimensions $n\gs 5$ for `efficient' illumination of convex bodies satisfying the assumptions of \cite[Lemma 8]{Tikhomirov-2017}, given that $2^{n-1}+2n < 2^n$ as soon as $n\gs 5$. However this is not the case either.

\begin{enumerate}
\item[2.] Consider the set $\B_2\in \mathcal{S}^9$ whose vertices are all the coordinate reflections and permutations of the points $e_1$ and $\frac{1}{2}(e_1+e_2+e_3+e_4)$. As previously, we can see that, if we want to illuminate all the coordinate reflections and permutations of $e_1$, we have to include all the directions $\pm e_i, i\in [9]$. On the other hand, among the outer normals at $\frac{1}{2}(e_1+e_2+e_3+e_4)$ are vectors of the form $e_j+e_k$ and of the form $e_j+e_k+e_s$ where $\{j,k,s\}$ is any subset of $[4]$ with cardinality 3. Thus, by Criterion A we can quickly check that none of the directions $\pm e_i$ illuminates $\frac{1}{2}(e_1+e_2+e_3+e_4)$, neither do any of the directions $\pm e_i\pm e_j$.

Moreover, we can check that $\frac{1}{2}(e_1+e_2+e_3+e_4)$ also has outer normals of the form 
\begin{gather*}
e_1+e_2+e_3+e_4\ \pm e_{t_1}\,\pm e_{t_2}\,\pm e_{t_3}\,\pm e_{t_4} 
\\
\hbox{where} \ \{t_1,t_2,t_3,t_4\}\subset [9]\setminus[4]\ \hbox{with cardinality 4}.
\end{gather*}
This shows that none of the directions in $T_1$ or $T_2\setminus\{\pm e_9\}$ (which all have at least 8 non-zero coordinates) can illuminate $\frac{1}{2}(e_1+e_2+e_3+e_4)$. We conclude that, for the coordinate reflections and permutations of $\frac{1}{2}(e_1+e_2+e_3+e_4)$ we need directions from $\{-1,0,1\}^9\setminus\{\vec{0}\}$ which have at least 3 and at most 7 non-zero coordinates.
\item[3.] Consider the set $\B_3\in \mathcal{S}^{25}$ whose vertices are all the coordinate reflections and permutations of the points
\begin{equation*}
e_1,\quad \frac{2}{3}(e_1+e_2+e_3),\quad \hbox{and}\ \frac{1}{3}\sum_{i=1}^{12}e_i.
\end{equation*}
Then we could illuminate $e_1$ (for example) using any of the directions $-e_1$, $-e_1\pm e_2$ or $-e_1\pm e_2\pm e_3$ (or even variations of the latter two): this is because e.g.
\begin{equation*}
e_1 + \frac{1}{2}(-e_1+e_2) = \frac{1}{2}(e_1+e_2)\in\intr\B,
\end{equation*}
as we can see if we compare coordinates with those of the vertex $\frac{2}{3}(e_1+e_2+e_3)$ (recall Lemma \ref{lem:smaller-coordinates}). We could also illuminate $\frac{2}{3}(e_1+e_2+e_3)$ using the directions $-e_1-e_2$ or $-e_1-e_2-e_3$ (whereas $-e_1$ would not work in this case). Moreover, among the outer normals at $\frac{2}{3}(e_1+e_2+e_3)$ we also have vectors of the form $e_1+e_2+e_3\pm e_{t_1}\pm e_{t_2}\pm e_{t_3}$ where $\{t_1,t_2,t_3\}$ is any subset of $[25]\setminus [3]$ with cardinality 3. This shows that any direction $d$ from $\{-1,0,1\}^{25}\setminus\{\vec{0}\}$ which illuminates $\frac{2}{3}(e_1+e_2+e_3)$ must have at most 5 non-zero coordinates, otherwise it would be orthogonal to one of the outer normals of $\frac{2}{3}(e_1+e_2+e_3)$, contradicting Criterion A.

On the other hand, vectors of the form $\sum_{s=1}^6 e_{j_s}$, where $\{j_s:1\ls s\ls 6\}$ is a subset of $[12]$ with cardinality 6, are among the outer normals at $\frac{1}{3}\sum_{i=1}^{12}e_i$, so any directions $d^\prime$ from $\{-1,0,1\}^{25}\setminus\{\vec{0}\}$ which illuminate this boundary point of $\B_3$ must have at least 7 and (as we can also check) at most 23 non-zero coordinates (so they cannot be among the directions from $T_1$ and $T_2$, nor can they be any of the directions we used for the second type of vertices).
\end{enumerate}

\medskip

The above examples indicate that we can construct `tricky' convex bodies, which cannot be covered by \cite[Lemma 8]{Tikhomirov-2017} (or even by \cite[Lemma 7]{Tikhomirov-2017}), in as high a dimension as we want. Furthermore, the set $T_2$ cannot be replaced by any set of the form
\begin{equation*}
\bigl(\{-1,1\}^{n-1}\times\{0\}\bigr) \cup R^\prime_0
\end{equation*}
where $R^\prime_0$ would only contain directions from $\{-1,0,1\}^n\setminus\{\vec{0}\}$ with a predetermined number of non-zero coordinates; or alternatively where $R^\prime_0$ would have cardinality $\ls p(n)$ with $p(x):\R\to \R$ being a fixed polynomial. On the other hand, as we explained in Subsection \ref{subsec:Tikh-and-our-method}, $R^\prime_0$ could be any set that has property (P1): e.g. Tikhomirov's probabilistic construction of a union of `good' realisations of the sets $\mathcal{S}_k, k=1,\ldots,\lceil\frac{n}{2}\rceil$, introduced on \cite[Page 378]{Tikhomirov-2017}, which is shown to have small enough cardinality for sufficiently high dimensions.

\medskip

There is one more point worth remarking: it can be checked that all the above, `tricky' examples satisfy ${\rm dist}(\B,[-1,1]^n) \gs 2$. This is no accident: if instead we had $1\neq {\rm dist}(\B,[-1,1]^n) < 2$, then one of the sets $T_1$ and $T_2$ ought to have worked, which is not the case for any of the above examples.

\medskip

{\bf Claim F.} Let $n\gs 3$ and let $\B\in \SU^n$ with $1\neq {\rm dist}(\B,[-1,1]^n) < 2$. Then $\B$ is illuminated by one of the sets $T_1$ and $T_2$ (using the notation in Proposition B). 
\smallskip \\
{\bf Note.} This is essentially the statement of \cite[Proposition 5]{Tikhomirov-2017} (that statement also includes the case $n=2$, which follows immediately from \cite[Lemma 8]{Tikhomirov-2017} for $n=2$, since $1\neq {\rm dist}(\B,[-1,1]^n)$ guarantees that $e_1+e_2\notin\B$ for $\B\in \SU^2$, and for such convex bodies the corresponding set $T_2 = \{\pm e_1,\pm e_2\}$ is easily seen to be an illuminating set).

\medskip

We defer the proof to Appendix A. Here we just mention one more example of a `tricky' convex body $\B\in \SU^n$ with $1\neq {\rm dist}(\B,[-1,1]^n) < 2$, which satisfies $\|e_1+e_2\|_\B>1$, and which nonetheless cannot be illuminated by the set $T_2$ (however, by Claim F we now know that it will be illuminated by $T_1$).

\begin{enumerate}
\item[4.] Let $n\gs 4$ and consider the set $\B_4^n\in \SU^n$ whose vertices are all the coordinate reflections and permutations of the points $e_1+\frac{3}{4}e_n$ and $\frac{1}{2}\sum_{i=1}^{n-1}e_i + \frac{3}{4}e_n$. Note that $\B_4$ contains the convex combination
\begin{equation*}
\frac{1}{n}\sum_{i=1}^n \left(\frac{3}{4}e_i \,+\! \sum_{j\in [n]\setminus\{i\}}\frac{1}{2}e_j\right) = \left(\frac{1}{n}\cdot\frac{3}{4} + \frac{n-1}{n}\cdot\frac{1}{2}\right)\sum_{i=1}^n e_i,
\end{equation*}
which shows that ${\rm dist}(\B,[-1,1]^n) \ls \frac{4n}{2n+1} < 2$.

Now observe that among the outer normals at $e_1+\frac{3}{4}e_n$ is the vector $e_1$, as well as vectors of the form $e_1+e_n\pm e_j$ where $j\in [n]\setminus\{1,n\}$. This shows that none of the directions in $T_2= \bigl(\{-1,1\}^{n-1}\times\{0\}\bigr)\cup\{\pm e_n\}$ can illuminate $e_1+\frac{3}{4}e_n$.
\end{enumerate}

\section{One illuminating set to `rule' them all}\label{sec:main-illuminating-set}

Here we will prove Theorem D, that is, we will show that we can find a single subset of the set 
\begin{equation*}
    G^n(\delta) := \left\{d\in \R^n: \exists\, i\in [n]\ \hbox{such that}\ d=\pm e_i\,+\!\sum_{j \in [n] \backslash \{i\}} 
   \! \pm \delta e_j \right\}
\end{equation*}
which is an `efficient' illuminating set for any 1-symmetric convex body in $\R^n$ (as long as $\delta$ is also chosen sufficiently small).
To do so, and motivate the methods by which we construct/select such subsets in Proposition \ref{prop:general-illuminating-set}, we single out the next two results and state them as separate lemmas; this will also help with the proofs in the next section.

\begin{lemma}\label{lem:simple-deep-illumination}
Let $n\gs 3$, $\B \in \SU^n $, and $x \in \partial \B$. Consider $k \in [n]$ such that $x_{k} \neq 0$. Then for any $\delta\in (0, \frac{1}{n})$, $x$ is illuminated by 
\begin{equation*}
   -\sign(x_{k})e_{k}\,+\!\! \sum_{i \in [n]\backslash \{k\}}\!\!\!\!(-\delta\sign(x_i))e_i
\end{equation*}
(here, if $x_i=0$, we simply make a choice for $\sign(x_i)$, setting it equal to either $+1$ or $-1$; no matter the choice of signs for these coordinates, the conclusion remains the same).
\end{lemma}
\begin{proof}
Denote the direction in the statement by $d$. WLOG, suppose that $k=1$ and also that $\abs{x_i} \geq \abs{x_{i+1}}$ for $2 \leq i \leq n-1$ (we can always reorder the indices in this argument if this is not the case; note also that the second assumption does not hinge on whether $\abs{x_1} = \norm{x}_\infty$ or not, we simply order the remaining coordinates of $x$).

If $\abs{\cZ_x}=0$ (that is, if $x$ has no zero coordinates), then the conclusion follows immediately by Corollary \ref{cor:uncond-illum}. 
Thus we now focus on the case where
\begin{equation*}
\abs{x_2}\gs \abs{x_3}\gs \cdots \gs \abs{x_{m-1}}\, >\,\abs{x_m}=0=\abs{x_{m+1}}=\cdots=\abs{x_n} 
\end{equation*}
for some $2\ls m\ls n$. Let us examine the entries of $x+|x_1|d$:
\begin{itemize}
\item $(x+|x_1|d)_1 = 0$,
\item for every $2\ls i\ls m-1$, $\abs{(x+|x_1|d)_i} = \max(|x_i| - |x_1|\delta, \,|x_1|\delta - |x_i|)$ (and if the maximum is not equal to the first argument, then $|x_i|$ must be `very small', $< \delta |x_1|$),
\item and for every $m\ls i\ls n$, $\abs{(x+|x_1|d)_i} = |x_1|\delta$.
\end{itemize}
Let $m^\prime$ be the smallest index $\gs 2$ such that $|x_{m^\prime}| < \frac{1}{n}|x_1|$.
We consider the following convex combination, which will be in $\B$ because of the 1-symmetry:
\begin{align*}
&\frac{1}{n-m^{\prime}+2}\Bigl[(|x_1|, |x_2|,|x_3|,\ldots,|x_{n-1}|,|x_n|) \,+\,  (|x_{m^{\prime}}|, |x_2|,|x_3|,\ldots, |x_{m^\prime-1}|,|x_1|,\ldots,|x_{n-1}|,|x_n|) 
\\
&\qquad + (|x_{m^{\prime}+1}|, |x_2|,|x_3|,\ldots, |x_{m^\prime-1}|,|x_{m^{\prime}}|, |x_1|,\ldots,|x_{n-1}|,|x_n|)\  +
\\
&\qquad\quad + \cdots\cdots + (|x_{n-1}|, |x_2|,|x_1|,\ldots,|x_1|,|x_n|) + (|x_n|, |x_2|,|x_1|,\ldots,|x_{n-1}|,|x_1|)\Bigr].
\end{align*}
We can observe the following regarding its entries:
\begin{itemize}
\item the 1st entry is $\frac{1}{n-m^{\prime}+2}\left(|x_1|+\sum_{i=m^{\prime}}^n |x_i|\right) \gs \frac{\abs{x_1}}{n-m^{\prime}+2} > 0$.
\item For every $2\ls i\ls m^{\prime}-1$, its $i$-th entry is equal to $|x_i|$.
\item For every $m^{\prime}\ls i\ls m-1$, its $i$-th entry is equal to 
\begin{equation*}
\frac{n-m^{\prime}+1}{n-m^{\prime}+2}|x_i| + \frac{1}{n-m^{\prime}+2}|x_1|\gs \frac{1}{n-m^{\prime}+2}|x_1|.
\end{equation*}
\item For every $m\ls i\ls n$, its $i$-th entry is equal to $\frac{1}{n-m^{\prime}+2}|x_1|$.
\end{itemize}
Thus, as long as we choose $\delta < \frac{1}{n}$, this convex combination will have strictly larger (in absolute value) respective entries compared to $x+|x_1|d$. 
By Lemma \ref{lem:smaller-coordinates} it follows that $x+|x_1|d\in \intr\B$.
\end{proof}

The following is a strengthening of the previous lemma which will be necessary in certain cases.

\begin{lemma}\label{lem:trickier-deep-illumination}
Let $n\gs 3$, $\B \in \SU^n $, and $x \in \partial \B$. Write $M_x$ for the set of maximum coordinates of $x$:
\begin{equation*}
M_x=\{k\ls n: \abs{x_k} = \norm{x}_\infty\},
\end{equation*}
and let $k_0\in M_x$. Then for any $\delta\in (0, \frac{1}{n})$ and for any \underline{arbitrary} choice of signs $\epsilon_i, i\notin M_x$, $x$ is illuminated by 
\begin{equation*}
   -\sign(x_{k_0})e_{k_0}\,+\!\! \sum_{k \in M_x\backslash \{k_0\}}\!\!\!\!(-\delta\sign(x_k))e_k \,+\!\!\sum_{i \in [n]\backslash M_x}\!\!\!\!\delta\epsilon_ie_i.
\end{equation*}
\end{lemma}
\begin{proof}
Fix a choice of signs $\epsilon_i, i\notin M_x$, and write $d$ for the corresponding direction. If $[n]\backslash M_x = \emptyset$, the conclusion follows immediately from Corollary \ref{cor:uncond-illum}, so assume there is some index $i\notin M_x$ and set $t=\abs{[n]\backslash M_x}$. For any such index $i$, let $\alpha_i := \abs{x_{k_0}} - \abs{x_i}$, and note that $\alpha_i \in (0,\norm{x}_\infty]$.
Set also
\begin{equation*}
\lambda_i : = \frac{\alpha_i^{-1}}{1+\sum_{j\in [n]\backslash M_x}\alpha_j^{-1}},\qquad \lambda_0 := \frac{1}{1+\sum_{j\in [n]\backslash M_x}\alpha_j^{-1}}.
\end{equation*}
Then $\lambda_0+\sum_{j\in [n]\backslash M_x}\lambda_j = 1$. Consider the following convex combination of points in $\B$:
\begin{align*}
w = \lambda_0\cdot \sum_{i=1}^n\abs{x_i}e_i + \sum_{j\in[n]\backslash M_x} \lambda_j\cdot \left(\abs{x_{k_0}}e_j + \abs{x_j}e_{k_0} + \sum_{i\notin \{k_0,j\}}\abs{x_i}e_i\right).
\end{align*}
Let us examine its coordinates:
\begin{itemize}
\item if $i\in M_x\backslash \{k_0\}$, then $w_i = \abs{x_i}$.
\item For $i=k_0$ we have $w_{k_0} = \lambda_0\abs{x_{k_0}} + \sum_{j\in[n]\backslash M_x}\lambda_j\abs{x_j}$.
\item If $i\in [n]\backslash M_x$, then $w_i = \lambda_i\abs{x_{k_0}} + (1-\lambda_i)\abs{x_i}$.
\end{itemize}
Thus for $i\notin M_x$, we have that $w_i > \abs{x_i}$, while $0 < w_{k_0} < \abs{x_{k_0}} = \norm{x}_\infty$. In particular,
\begin{equation*}
w_i - \abs{x_i} = \lambda_i\cdot \bigl(\abs{x_{k_0}} - \abs{x_i}\bigr) = \lambda_i\cdot \alpha_i = \frac{1}{1+\sum_{j\in [n]\backslash M_x}\alpha_j^{-1}} = \lambda_0,
\end{equation*}
while
\begin{equation*}
\abs{x_{k_0}} -w_{k_0} = \sum_{j\in[n]\backslash M_x}\lambda_j\cdot\bigl(\abs{x_{k_0}}-\abs{x_j}\bigr) = \sum_{j\in[n]\backslash M_x}\lambda_j\cdot\alpha_j = t\cdot\lambda_0.
\end{equation*}
From this it's clear that, for each $i\notin M_x$, 
\begin{equation}\label{eqp1:lem:trickier-deep-illumination}
w_i - \abs{x_i} = \frac{1}{t}\bigl(\abs{x_{k_0}} -w_{k_0}\bigr) > \frac{1}{n}\bigr(\abs{x_{k_0}} -w_{k_0}\bigr).
\end{equation}

Consider now the vector $y=x+\bigl(\abs{x_{k_0}} -w_{k_0}\bigr)d$. Let us examine its coordinates too:
\begin{itemize}
\item $y_{k_0} = x_{k_0} - \sign(x_{k_0})\bigl(\abs{x_{k_0}} -w_{k_0}\bigr) = \sign(x_{k_0})w_{k_0}$.
\item If $i\in M_x\backslash \{k_0\}$, then 
\begin{equation*}
y_i = x_i -  \delta\sign(x_i)\bigl(\abs{x_{k_0}} -w_{k_0}\bigr) = \sign(x_i)\bigl(\abs{x_i}-\delta\abs{x_{k_0}} +\delta w_{k_0}\bigr)
= \sign(x_i)\big((1-\delta)\abs{x_{k_0}} + \delta w_{k_0}\bigr),
\end{equation*}
which is smaller in absolute value than $\abs{x_{k_0}} = \abs{x_i} = w_i$.
\item Finally, if $i\in [n]\backslash M_x$,
\begin{equation*}
\abs{y_i} \ls \abs{x_i} + \delta\bigl(\abs{x_{k_0}} -w_{k_0}\bigr).
\end{equation*}
Given that $\delta < \frac{1}{n}$, by \eqref{eqp1:lem:trickier-deep-illumination} we obtain that
\begin{equation*}
\abs{x_i} + \delta\bigl(\abs{x_{k_0}} -w_{k_0}\bigr) < \abs{x_i} + \frac{1}{n}\bigl(\abs{x_{k_0}} -w_{k_0}\bigr) < w_i.
\end{equation*}
\end{itemize}
Thus all coordinates except the $k_0$-th one of $y$ are strictly smaller in absolute value than the respective coordinates of $w$ (while $\abs{y_{k_0}} = w_{k_0} < 1$). It remains to consider the section
\begin{equation*}
\bigl\{\xi\in \B: \xi_{k_0} = y_{k_0}\bigr\} = \bigl\{\xi\in \B: \xi_{k_0} = \sign(y_{k_0})w_{k_0}\bigr\}
\end{equation*}
and apply Corollary \ref{cor:affine-set}.
\end{proof}

\begin{remark}\label{rem:generalise-deep-illumination}
If we look again at the proof of Lemma \ref{lem:simple-deep-illumination}, then we can also derive the following more general conclusion: if $x\in \partial \B$ and $x_k\neq 0$, then $x$ is illuminated by any direction of the form
\begin{equation*}
   -\sign(x_{k})e_{k}\,+\!\! \sum_{i \in [n]\backslash \{k\}}\!\!\!\!(-\delta_i\sign(x_i))e_i
\end{equation*}
where $\sign(x_i)$ is arbitrarily chosen from $\{-1,1\}$ for $i\in \cZ_x$ and where $\delta_i\in (0,\frac{1}{n})$ for all $i\in [n]\setminus\{k\}$ (namely we don't have to set all the `small' coordinates equal in absolute value, we just have to make sure all of them are $< \frac{1}{n}$).

\medskip

Analogously, by examining the proof of Lemma \ref{lem:trickier-deep-illumination}, we can see that, if $k_0\in M_x$, then $x$ can be illuminated by any direction of the form 
\begin{equation*}
   -\sign(x_{k_0})e_{k_0}\,+\!\! \sum_{k \in M_x\backslash \{k_0\}}\!\!\!\!(-\delta_k\sign(x_k))e_k \,+\!\!\sum_{i \in [n]\backslash M_x}\!\!\!\!\widetilde{\delta}_i\,\epsilon_ie_i
\end{equation*}
where $\epsilon_i\in \{-1,1\}$ is arbitrary for $i\in [n]\backslash M_x$, and where $\delta_k,\widetilde{\delta}_i \in (0,\frac{1}{n})$ (but not necessarily the same number). In fact, for $i\in [n]\backslash M_x$, we can even allow $\widetilde{\delta}_i$ to be chosen from $[0,\frac{1}{n})$, namely allow $\widetilde{\delta}_i$ to be set equal to 0.

Note that this more general version of Lemma \ref{lem:trickier-deep-illumination} extends Lemma \ref{lem:unique-max-coord}.
\end{remark}

\medskip

\begin{proposition}\label{prop:general-illuminating-set}
Let $n \gs 2$, $\delta \in (0,1)$, and consider the set $G^n=G^n(\delta)$ from \eqref{def:big-illum-set} (or from the first paragraph of Section \ref{sec:main-illuminating-set}).

There exists a subset $\I^n \subset G^n$ with $\abs{\I^n}=2^n$ and with the property: for any vector $x \in \R^n\backslash\{\vec{0}\}$, there is $d \in \I^n$ such that $\sign(x_i) = -\sign(d_i)$ for every $i \in [n]\backslash\cZ_x$ and moreover $1=\norm{d}_\infty =\abs{d_j}$ for an index $j \in [n]\backslash\cZ_x$.
\end{proposition}

\begin{definition} We will be using the following terminology in settings such as the above: given a vector $x \in \R^n\backslash\{\vec{0}\}$ and given a set $\I \subset G^n$, we will say that $\I$ \textit{deep illuminates $x$} if there exists $d \in \I$ such that $\sign(x_i) = -\sign(d_i)$ for every $i \in [n]\backslash\cZ_x$ and at the same time $1=\norm{d}_\infty = \abs{d_j}$ for an index $j \in [n]\backslash \cZ_x$ (moreover, if $\I=\{d\}$ is a singleton, we will say that the direction $d$ deep illuminates $x$).
\end{definition}
\begin{proof} We will present two methods of finding sets $\I^n$ with the desired properties, with the first method being more geometric, while the second one will be combinatorial. In most of the results in the sequel, we can use either method to produce a set $\I^n$ with the desired properties. Note however that only the 2nd method here 
\begin{itemize}
\item[--] ensures that the set $\I^n$ that we get is also formed from pairs of opposite directions, 
\item[--] and moreover it provides us with a better overview of which directions from $G^n$ end up inside $\I^n$. 
\end{itemize}
These will be important in some of our final results.

\medskip

 {\bf Method 1.} Consider the graph that we get from the vertices and edges of $[-1,1]^n$. Recall that we can find Hamilton cycles on this graph, that is, closed and simple paths that pass by all the vertices of $[-1,1]^n$ (exactly once since they are simple) and return to their initial vertex. More simply, a Hamilton cycle of the graph of the cube is equivalent to a sequence/ordering $w_1, w_2,..., w_{2^n}$ of the vertices
\begin{itemize}
\item such that all vertices appear in the sequence 
\item and such that consecutive vertices in the sequence differ in exactly one coordinate; this should be true for the pair of vertices $w_{2^n}$ and $w_1$ as well.
\end{itemize}
Note also that, since the path is closed, we could start at any vertex. Furthermore, applying a symmetry of the cube to such a path (that is, any composition of coordinate permutations and of sign changes) will just give us another path with the same properties. Thus it is fine to require that $w_1 = (1,1,1,..., 1,1)$, while $w_{2^n} = (1,1,1,..., 1,-1)$.

To give an example of such a Hamilton cycle/ordering of the vertices of $[-1,1]^n$, let us start from the sequence
\begin{align}
\label{D2 Sequence}
    w_1=(1,1), \ w_2 =(-1,1), \ w_3=(-1,-1), \ w_4=(1,-1),
\end{align}
which satisfies the properties we want when $n=2$. We can construct admissible sequences in higher dimensions via recursion.
Given $n> 2$, assume that we have already found/fixed such a sequence in dimension $n-1$, and that we denote it by $P^{n-1} = \bigl(w_1^{n-1}, w_2^{n-1},..., w_{2^{n-1}}^{n-1}\bigr)$. As explained before, WLOG we can assume that 
\begin{equation*}
w_1^{n-1} = \sum_{i=1}^{n-1} e_i\,, \qquad \hbox{while}\quad w_{2^{n-1}}^{n-1} = - e_{n-1} + \sum_{i=1}^{n-2}e_i
\end{equation*}
(the $e_i$'s here are the standard basis vectors in $\R^{n-1}$).
Then an admissible sequence for dimension $n$ is the following:
\begin{align*}
&w_1^n = \bigl(w_1^{n-1}, 1\bigr), \ \ w_2^n = \bigl(w_2^{n-1}, 1\bigr),\  \ldots\\
&\hspace{2cm}\ldots\,, \ w_{2^{n-1}}^n = \bigl(w_{2^{n-1}}^{n-1}, 1\bigr),\ \ w_{2^{n-1}+1}^n = \bigl(w_{2^{n-1}}^{n-1}, -1\bigr),
\\
&w_{2^{n-1}+2}^n = \bigl(w_{2^{n-1}-1}^{n-1}, \,-1\bigr),\ \ w_{2^{n-1}+3}^n = \bigl(w_{2^{n-1}-2}^{n-1}, \,-1\bigr), \ \ldots\\
&\hspace{2cm}\ldots\,,\  w_{2^n-1}^n = \bigl(w_2^{n-1},\,-1\bigr), \ \ w_{2^n}^n = \bigl(w_1^{n-1},\,-1\bigr).
\end{align*}
In other words, we append to every vector in the sequence $P^{n-1}$ one more coordinate equal to $1$, and then we revert the order of the terms of $P^{n-1}$, and to each of the repeated and reverse-reordered terms we append one more coordinate equal to $-1$. Finally we concatenate the two new subsequences that we got to end up with a sequence $P^n$ formed by $2^n$ vectors in $\R^n$ (all of which are vertices of $[-1,1]^n$).

It is not hard to check that the sequence $P^n$ has the desired properties (assuming that the sequence $P^{n-1}$, which we started with, did so as well).

\smallskip

\noindent {\bf Note.} Such orderings of the vertices of $[-1,1]^n$ are sometimes also called \emph{Gray codes}, and they play an important role in digital communications and in analog-to-digital signal conversion. The examples that we gave in particular, which are constructed by this recursive procedure, are often called the \emph{binary-reflected Gray codes}.

\bigskip

Next we explain how we construct a desirable subset $\I^n$ of $G^n$ once we have a sequence $P^n = (w_1^n, w_2^n, ..., w_{2^n}^n)$ as above. Fix an index $r\in \{1,2,..., 2^n\}$), and consider the \underline{ordered} pair of vertices $(w_r, w_{r+1})$  (from now on we will be suppressing the superscript $n$ which indicates the dimension). Note that, if $r=2^n$, then the pair to consider is $(w_{2^n}, w_1)$. This leads to the direction 
\begin{equation*}
\mathbf{d}_r : = -w_{r,i_r}e_{i_r} + \sum_{i\in [n]\setminus \{i_r\}} (-\delta\cdot w_{r,i})e_i, 
\end{equation*}
where $i_r$ is the index of the unique entry in which $w_r$ and $w_{r+1}$ differ (in fact, $-w_{r,i_r} = w_{r+1,i_r}$).
For instance, $\mathbf{d}_{2^n} = (-\delta, -\delta,..., -\delta,-\delta, +1)$. See also a visualisation of these directions in $\R^3$.
\begin{center}
        \includegraphics[width=5cm]{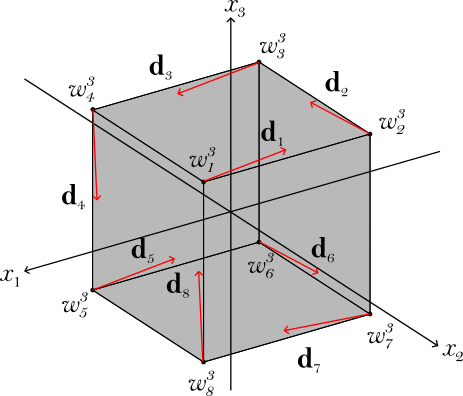}
        \\
        {\bf \footnotesize Figure 1:} \emph{\footnotesize A Hamilton cycle on the graph of $[-1,1]^3$,} 
        \\
        \emph{\footnotesize and the corresponding directions $\mathbf{d}_1, \mathbf{d}_2,\ldots, \mathbf{d}_8$ for some $\delta \ll 1$.}
    \end{center}
    
    \smallskip

We show that the set
\begin{equation*}
\I^n \equiv \I^n(\delta) : = \{\mathbf{d}_1, \mathbf{d}_2,\ldots, \mathbf{d}_{2^n -1},\mathbf{d}_{2^n}\}
\end{equation*}
deep illuminates any vector $x\in {\mathbb R}^n\setminus\{\vec{0}\}$. Fix such a vector $x$, and consider the vector $y$ of its coordinate signs:
\begin{equation*}
y\equiv y_x : = \bigl(\widetilde{\sign}(x_1),\widetilde{\sign}(x_2), \ldots, \widetilde{\sign}(x_{2^n-1}), \widetilde{\sign}(x_{2^n})\bigr),
\end{equation*}
where $\widetilde{\sign}(x_i) = \sign(x_i)$ if $x_i\neq 0$, and $\widetilde{\sign}(x_i) = 0$ if $x_i=0$. 

Clearly $y_x$ is found on the boundary of the cube $[-1,1]^n$, and if $x$ has no zero coordinates, then $y_x$ is actually a vertex of the cube. Otherwise there is a unique face $F_y$ of the cube with dimension $k\in \{1,2,...,n-1\}$ such that $y$ is in the relative interior of $F_y$.

Consider now the vertices of the cube which belong to $F_y$ and call them `good' for $y$. More simply, these are precisely the vertices $v_j$ of the cube which agree with $y$ in all entries where $y$ has a non-zero coordinate: $v_j\in F_y$ if and only if $v_{j,i} = y_i$ for all $i\in [n]\setminus \cZ_y = [n]\setminus\cZ_x$. Observe that there will be $2^k$ vertices which are `good' for $y$, where $k$ is the number of zero coordinates of $y$, $k=\abs{\cZ_y} = \abs{\cZ_x}$. Since $x\in {\mathbb R}^n\setminus\{\vec{0}\}$, there will be at most $2^{n-1}$ `good' vertices for $y$.

{\bf Claim.} There exists an index $r_0\in \{1,2,..., 2^n\}$ such that the $r_0$-th vertex $w_{r_0}$ in the sequence $P^n$ is `good' for $y$, while the next vertex $w_{r_0+1}$ is `bad' for $y$ (here we take $r_0+1 \ \hbox{mod}\ 2^n$, that is, $r_0+1=1$ if $r_0=2^n$).

\emph{Proof of the claim.} Since there exists at least one and at most $2^{n-1}$ `good' vertices for $y$, we can set $s_0\in\{1,2,...,2^n\}$ to be the largest index of a `good' vertex for $y$, and we can also set $t_0$ to be the smallest index of a `bad' vertex for $z$. We distinguish the following cases.
\begin{itemize}
\item[1.] $s_0 < 2^n$. Then $s_0 + 1\ls 2^n$ and is larger than $s_0$, so it cannot be the index of a `good' vertex for $y$. Thus we can set $r_0=s_0$, since $w_{s_0}$ is a `good' vertex for $y$ and $w_{s_0+1}$ is a `bad' vertex for $y$.
\item[2a.] $s_0=2^n$ and $t_0=1$. By our choices for $s_0$ and $t_0$, we have that $w_{s_0}$ is a `good' vertex for $y$ and $w_{t_0}$ is a `bad' vertex for $y$. Moreover, in this case $s_0+1 \ \hbox{mod}\ 2^n = t_0 \ \hbox{mod}\ 2^n$, so we can set $r_0=s_0$ again.
\item[2b.] $t_0>1$. Then $t_0-1\gs 1$ and is smaller than $t_0$, so, by our choice for $t_0$, $t_0-1$ cannot be the index of a `bad' vertex for $y$ (while $w_{t_0}$ is indeed a `bad' vertex for $y$). Thus we can set $r_0=t_0-1$. 
\end{itemize}

Finally, we verify that, if $r_0\in \{1,2,..., 2^n\}$ is an index which satisfies the property in the claim, then the direction $\mathbf{d}_{r_0}$ deep illuminates the points $y$ and $x$. Indeed, recall that respective entries of $\mathbf{d}_{r_0}$ and of the vertex $w_{r_0}$ have opposite signs. Since $w_{r_0}$ is `good' for $y$, for every $i\in [n]\setminus\cZ_y$ we have that $y_i = w_{r_0,i} = -\sign(\mathbf{d}_{r_0,i})$. On the other hand, $w_{r_0+1}$ is a `bad' vertex for $y$, so we can find some index $i_0$ such that $y_{i_0}\neq 0$ and $y_{i_0}\neq w_{r_0+1,i_0}$. But $w_{r_0}$ and $w_{r_0+1}$ are consecutive vertices in the sequence $P^n$, so they can only differ in exactly one entry; then this must be the $i_0$-th entry. By construction of the directions $\mathbf{d}_r$, we see that $\norm{\mathbf{d}_{r_0}}_\infty$ is attained in the $i_0$-th entry for which we have $y_{i_0}\neq 0$. These combined show that $\mathbf{d}_{r_0}$ deep illuminates $y$ (which is equivalent to saying that $\mathbf{d}_{r_0}$ deep illuminates $x$).

This completes the proof of the proposition by the 1st method.

\bigskip

{\bf Method 2.} We now discuss a more direct, and fully combinatorial, way of constructing of $\I^n=\I^n(\delta)$ (which will also rely on a recursive definition). Using this method, we will see that we can also ensure that 
\begin{center}
the resulting set $\I^n$ consists of pairs of opposite directions
\end{center}
(compare with Figure 1, and the set $\I^3(\delta)$ that the 1st method would give us, which does not have this property).

We start by setting
\begin{equation*}
\I^2 = \{(+1,+\delta), (-1,-\delta), (+\delta,-1), (-\delta,+1)\} = \bigl\{\pm (1,\delta),\,\pm (\delta, -1)\bigr\}.
\end{equation*}

Next, out of the $8=2^3$ directions that $\I^3$ will eventually contain, we define the first 4 to be the following:
\begin{equation*}
((+1,+\delta),+\delta),\  ((-1,-\delta),-\delta),\  ((+\delta,-1),-\delta),\  ((-\delta,+1),+\delta).
\end{equation*}
In other words, we form these directions by appending one more `small' coordinate to the vectors in $\I^2$ in such a way that the $3$rd coordinate will have the same sign as the $2$nd coordinate. 

\smallskip

\noindent {\bf Important Note.} This also causes the last coordinate of the new direction to have the same sign as the maximum (in absolute value) coordinate of the preexisting direction (we should keep this in mind as it will help us with the recursive definition and also later, when checking the desired properties of the sets we define). 

\smallskip

Note finally that, since $\I^2$ consisted of two pairs of opposite directions, these new directions also form pairs of opposite directions.

\medskip

The remaining 4 directions of $\I^3$ are defined as follows: in all these directions the maximum (in absolute value) coordinate will be the last one; moreover each of these directions will be formed 
\begin{itemize}
\item[--] by considering each of the first 4 directions at a time, 
\item[--] and by keeping the signs of its coordinates, except for the last sign which will be flipped (and thus will now be opposite to the sign of the second to last coordinate). 
\end{itemize}
In other words, the last 4 directions of $\I^3$ are set equal to
\begin{equation*}
(+\delta,+\delta,-1),\  (-\delta,-\delta,+1),\  (+\delta,-\delta,+1),\  (-\delta,+\delta,-1).
\end{equation*}
Again these form two pairs of opposite directions. 

\medskip

To illustrate the construction with one more example: the first 8 directions of $\I^4$ will be
\begin{gather*}
((+1,+\delta,+\delta),+\delta),\  ((-1,-\delta,-\delta),-\delta),\  ((+\delta,-1,-\delta),-\delta),\  ((-\delta,+1,+\delta),+\delta),
\\
\qquad\quad ((+\delta,+\delta,-1),-\delta),\  ((-\delta,-\delta,+1),+\delta),\  ((+\delta,-\delta,+1),+\delta),\  ((-\delta,+\delta,-1),-\delta),
\end{gather*}
while the following 8 will be
\begin{gather*}
((+\delta,+\delta,+\delta),-1),\  ((-\delta,-\delta,-\delta),+1),\  ((+\delta,-\delta,-\delta),+1),\  ((-\delta,+\delta,+\delta),-1),
\\
\qquad\quad ((+\delta,+\delta,-\delta),+1),\  ((-\delta,-\delta,+\delta),-1),\  ((+\delta,-\delta,+\delta),-1),\  ((-\delta,+\delta,-\delta),+1).
\end{gather*}
Again, we observe that the important features here are: 
\begin{itemize}
\item[(i)] the last coordinate is either the maximum in absolute value, or alternatively it has the same sign as the maximum coordinate;
\item[(ii)] the projections of the first 8 directions onto their first $(n-1)=3$ coordinates form precisely $\I^3$, while the corresponding projections of the last 8 directions have the same coordinate signs as the directions in $\I^3$ (and thus all combinations of three signs appear);
\item[(iii)] in the first 8 directions the last coordinate has the same sign as the second to last, while in the last 8 directions they have opposite signs;
\item[(iv)] the $16=8+8$ directions form 8 pairs of opposite directions.
\end{itemize}

Assuming now that $\I^{n-1}(\delta) = \I^{n-1}$ has been defined in a way that guarantees the analogous properties, we can form $\I^n(\delta)=\I^n$ as follows: out of the $2^n$ directions that $\I^n$ needs to have, the first $2^{n-1}$ will be formed by appending to each direction $\mathbf{d}^{n-1}_s$ of $\I^{n-1}$ one more `small' coordinate at the end, so that this new coordinate will have the same sign as the last coordinate of $\mathbf{d}^{n-1}_s$; that is,
\begin{equation*}
\mathbf{d}^n_s : = \bigl(\mathbf{d}^{n-1}_s,\,\sign(\mathbf{d}^{n-1}_{s,n-1})\delta\bigr).
\end{equation*} 
At the same time, $\mathbf{d}^{n-1}_s$ enables us to also define one of the remaining $2^{n-1}$ directions for $\I^n$, which we will denote by $\mathbf{d}^n_{2^{n-1}+s}$: the sign of each of the first $n-1$ coordinates of $\mathbf{d}^n_{2^{n-1}+s}$ will be the same as for the respective coordinate of $\mathbf{d}^{n-1}_s$, while the last coordinate of $\mathbf{d}^n_{2^{n-1}+s}$ will be equal to 1 in absolute value and will have opposite sign to the previous coordinate, the $(n-1)$-th one. That is,
\begin{equation*}
\mathbf{d}^n_{2^{n-1}+s} := \Bigl(\delta\cdot\Bigl(\sign(\mathbf{d}^{n-1}_{s,1}),\, \sign(\mathbf{d}^{n-1}_{s,2}),\ldots,\sign(\mathbf{d}^{n-1}_{s,n-2}),\,\sign(\mathbf{d}^{n-1}_{s,n-1})\Bigr),\ \,-\!\!\!\!-\sign(\mathbf{d}^{n-1}_{s,n-1})\Bigr).
\end{equation*}
Based on this definition and the assumption that $\I^{n-1}$ already has properties analogous to the ones we stated for $\I^4$, it is not hard to check that $\I^n$ continues to have the analogous properties.

\bigskip

It remains to verify that $\I^n$ satisfies the statement of Proposition \ref{prop:general-illuminating-set}, that is, that it deep illuminates every vector $x\in \R^n\backslash \{\vec{0}\}$. Just as with the above properties, we will verify this inductively.

Clearly $\I^2$ deep illuminates $\R^2$. Assume that, for some $n-1\gs 2$, $\I^{n-1}$ has already been defined 
in some way that ensures that it deep illuminates $\R^{n-1}$. Moreover, assume that $\I^n$ is defined with the help of this $\I^{n-1}$ as described in the last paragraph \emph{(note that we don't even need to know that $\I^{n-1}$ has the analogous properties (i)-(iv), simply that it deep illuminates $\R^{n-1}$, and that $\I^n$ is defined through it as described last)}. Consider $x\in \R^n\backslash \{\vec{0}\}$.
\begin{itemize}
\item If $x_n$ is the only non-zero coordinate of $x$, then we can use (half of) the last $2^{n-1}$ directions of $\I^n$ to deep illuminate $x$ (based on what $\sign(x_n)$ is).
\item If $x_n=0$ instead, then ${\rm Proj}_{e_n^\perp}(x)\in \R^{n-1}\backslash\{\vec{0}\}$ (we slightly abuse notation here), and thus there is $\mathbf{d}^{n-1}_s\in \I^{n-1}$ which deep illuminates this projection. But then $\mathbf{d}^n_s = \bigl(\mathbf{d}^{n-1}_s,\,\sign(\mathbf{d}^{n-1}_{s,n-1})\delta\bigr)$ deep illuminates $x$.
\item Finally, if both $x_n\neq 0$ and ${\rm Proj}_{e_n^\perp}(x)\in \R^{n-1}\backslash\{\vec{0}\}$, then again we find $\mathbf{d}^{n-1}_s\in \I^{n-1}$ which deep illuminates the projection. It follows that either $\mathbf{d}^n_s = \bigl(\mathbf{d}^{n-1}_s,\,\sign(\mathbf{d}^{n-1}_{s,n-1})\delta\bigr)$ or $\mathbf{d}^n_{2^{n-1}+s}$ will deep illuminate $x$ (since the only difference in coordinate signs between these two directions occurs in the last coordinate).
\end{itemize}
This completes the 2nd proof to Proposition \ref{prop:general-illuminating-set} too. \end{proof}

\medskip

We can now establish Theorem D from Subsection \ref{subsec:Tikh-and-our-method}.

\begin{theorem}\label{thm:illum-everything}
Let $n\gs 3$ and let $\I^n(1/(n+1)) \subset G^n(1/(n+1))$ be any of the sets satisfying the statement of Proposition \ref{prop:general-illuminating-set}. Then $\I^n(1/(n+1))$ illuminates all $1$-symmetric convex bodies in $\R^n$.
\end{theorem}
\begin{proof} Consider an arbitrary $1$-symmetric convex body $\B$ in $\R^n$.
Given our terminology, any set $\I^n(\delta)$ satisfying the statement of Proposition \ref{prop:general-illuminating-set} will deep illuminate all the non-zero vectors of $\R^n$, and hence all the boundary points of $\B$. Since we assume here that $\delta= \frac{1}{n+1} < \frac{1}{n}$, Lemma \ref{lem:simple-deep-illumination} finishes the proof.
\end{proof}

\begin{remark}\label{rem:thm:illum-everything}
To also cover the case $n=2$, we check directly that $\I^2(\delta_0) = \bigl\{\pm(1,\delta_0),\,\pm(-\delta_0, 1)\bigr\}$ illuminates any 1-symmetric convex body $\B$ in $\R^2$, as long as $\delta_0\in (0,1)$. 

Indeed, WLOG we can assume that $\B\in \SU^2$. If $x\in \partial \B$ has only one non-zero coordinate, say $x_1\neq 0$, then necessarily $x=\pm e_1$, and we will have that $\pm e_1 + \mp(1,\delta_0) = \mp \delta_0 e_2 \in \intr\B$, since $\delta_0 < 1$. Similarly we deal with the case where $x_2\neq 0$. 

On the other hand, if both coordinates of $x$ are non-zero, then we simply pick the unique direction from $\I^2(\delta_0)$ which has opposite corresponding signs to the signs of $x$; by Corollary \ref{cor:uncond-illum} this illuminates $x$.
\end{remark}

\begin{remark}\label{rem:deep-2-classical}
We can summarise the results of this section in the following statement: in the context of 1-symmetric convex bodies in $\R^n$, and of subsets (or elements) of $G^n(\delta)$, where $\delta\in (0,\frac{1}{n})$, we have that 
\begin{center}\emph{if a subset (or direction) deep illuminates, then it also illuminates.}\end{center}
\end{remark}

\section{Main illumination results}\label{sec:main-results}

Let $n\gs 2$ and let $\B$ be a 1-symmetric convex body in $\R^n$. Recall that we have defined $m_\B$ to be the largest $k\in[n]$ such that $\norm{e_1+e_2+\cdots+e_k}_\B = \norm{e_1}_\B$, and also that
\begin{equation*}
{\rm dist}(\B,\,[-1,1]^n) = \frac{\norm{e_1+e_2+\cdots+e_n}_\B}{\norm{e_1}_\B}
\end{equation*}
(following from the same definition of ${\rm dist}(\cdot, [-1,1]^n)$ as the one given in \cite{Tikhomirov-2017}). Thus, if $m_\B < n$, we have that ${\rm dist}(\B,\,[-1,1]^n) > 1$. Unfortunately this does not necessarily imply that one can hope to illuminate $\B$ with fewer than $2^n$ directions, because e.g. the unit cross-polytope $CP_1^2$ in $\R^2$, that is, the body ${\rm conv}\{\pm(1,0), \pm (0,1)\}$, satisfies ${\rm dist}(CP_1^2,\,[-1,1]^2) = 2$, but it is just a rotated square, so we need 4 directions to illuminate it.

That said, such examples cannot arise in higher dimensions {\bf in the class of 1-symmetric convex bodies.} Indeed, if $n\gs 3$ and $\B$ is a 1-symmetric body in $\R^n$ with ${\rm dist}(\B,\,[-1,1]^n) > 1$ (equivalently, such that $m_\B < n$), then, as we will shortly see, $\II(B) \ls 2^n - 2$.

\smallskip

For simplicity we also work with the normalisation $\B\in \SU^n$, that is, $\norm{e_i}_\B =1$ for all $i\in [n]$. We first prove the counterparts to \cite[Lemma 8 and Proposition 10]{Tikhomirov-2017}. Recall that these concern bodies $\B\in \SU^n$ which do not contain any `large' unit subcubes.

\subsection{Special cases: no large unit subcubes}

\begin{theorem}\label{thm:illum-no-unit-squares}
Let $n\gs 3$, and suppose $\B\in{\mathcal S}^n$ is such that $\norm{e_i+e_j}_\B >1$ for all $1\ls i < j\ls n$ (in other words, $m_\B =1$). Then $\B$ can be illuminated by the set 
\begin{equation*}
\bigl[\I^{n-1}\bigl(\tfrac{1}{n}\bigr)\times\{0\}\bigr]\cup \{\pm e_n\},
\end{equation*}
where $\I^{n-1}(1/n)$ is any of the sets given by Proposition \ref{prop:general-illuminating-set}. Thus $\I(\B) \ls 2^{n-1} + 2$.
\end{theorem}
\begin{proof}
Let $x\in \partial \B$. We distinguish two cases.
\begin{itemize}
\item If $|x_n| < 1$, then we could use Corollary \ref{cor:affine-set} with the affine set $\{\xi\in \R^n: \xi_n = x_n\}$ which contains the point $x_n e_n\in \intr\B$. Note that the projection of the section $\B_{x_n} = \B \cap \{\xi\in \R^n: \xi_n = x_n\}$ onto the first $n-1$ coordinates is a 1-symmetric convex body in $\R^{n-1}$ (allowing some abuse of notation), and ${\rm Proj}_{e_n^\perp}(x)$ is on its boundary (because if it weren't, $x$ would be in $\intr\B$ by Lemma \ref{lem:affine-set}). By Theorem \ref{thm:illum-everything} (and Remark \ref{rem:thm:illum-everything} in the case that $n=3$), we know that ${\rm Proj}(\B_{x_n})$ is illuminated by $\I^{n-1}(1/n)$. This implies that we can find $\mathbf{d}^{n-1}_s\in \I^{n-1}(1/n)$ such that 
\begin{equation*}
{\rm Proj}_{e_n^\perp}(x) + \varepsilon \,\mathbf{d}^{n-1}_s \in \intr{{\rm Proj}(\B_{x_n})}
\end{equation*}
for some $\varepsilon > 0$. Equivalently $x+ \varepsilon \,(\mathbf{d}^{n-1}_s,0)$ is in ${\rm relint}\, \B_{x_n}$, and thus in $\intr{\B}$.
\item If $|x_n|=1$, then, by the assumption that $\norm{e_i+e_j}_\B >1$ for all $1\ls i < j\ls n$, we have that $|x_i| < 1 = |x_n|$ for all $i\ls n-1$. But then we can use Lemma \ref{lem:unique-max-coord} (this case is handled in the same way as in \cite[Lemma 8]{Tikhomirov-2017}).
\end{itemize} 
The proof is complete.
\end{proof}

\medskip

Next we provide an alternative to \cite[Proposition 10]{Tikhomirov-2017}, which dealt with the illumination of those 1-symmetric convex bodies $\B$ that have `large' distance to the cube: ${\rm dist}(\B,\,[-1,1]^n) \gs 2$. Recall that this implies that $m_\B \ls \frac{n}{2}$. Now, using this last assumption as our starting point, we establish an alternative result which works for all dimensions $n\gs 4$ (observe that, if $n=3$, then $m_\B\ls \frac{n}{2}$ implies that $m_\B = 1$, thus the case $n=3$ is already handled by Theorem \ref{thm:illum-no-unit-squares}).

We first gather some estimates for certain partial sums of binomial coefficients (these will come up as upper bounds on the cardinalities of the illuminating sets we will use).

\smallskip

\begin{lemma}\label{lem:subset-cardinality}
For all $n\gs 3$, we have that
\begin{equation*}
2\sum_{s=0}^{\lfloor (n-1)/4\rfloor}\binom{n}{s} \ls 2^{n-1}-2.
\end{equation*}
Moreover, this inequality is sharp only when $n=3$. For higher $n$, we can obtain the following: if $n=5$, 
then the left-hand side is $=2^{n-1}-4$, if $n=6$, the left-hand side is $= 2^{n-1}-18$, while for all other $n$ it holds that
\begin{equation*}
2\sum_{s=0}^{\lfloor (n-1)/4\rfloor}\binom{n}{s} \ls 2^{n-1}-\binom{n}{\lfloor \tfrac{n}{2}\rfloor}.
\end{equation*}
Finally, for sufficiently high dimensions 
the left-hand side is even smaller: we have that
\begin{equation*}
2\sum_{s=0}^{\lfloor (n-1)/4\rfloor}\binom{n}{s} = \left(\frac{256}{27}\right)^{\frac{n}{4}\,+\,O(\log n)}.
\end{equation*}
\end{lemma}
\begin{proof}
For $n$ between 3 and 8, we can check the relevant inequalities directly. When $n\gs 9$, for convenience we will consider four cases:
\begin{description}
\item[$\bm{n=4k}$ for some $\bm{k\gs 3}$.] Then $\lfloor\frac{n-1}{4} \rfloor = k-1\gs 2$. We have that
\begin{equation*}
2^n = 2 \sum_{s=0}^{2k-1}\binom{n}{s} \ +\ \binom{n}{2k} \ = \ 2\sum_{s=0}^{k-1}\binom{n}{s} + 2\sum_{s=k}^{2k-1}\binom{n}{s} \ +\ \binom{n}{2k}. 
\end{equation*}
We have that $2\binom{n}{2k-1} = 2\frac{n!}{(2k-1)!(2k+1)!} = 2\frac{2k}{2k+1}\frac{n!}{((2k)!)^2} > \binom{n}{2k}$, and moreover we have that
\begin{equation*}
\binom{n}{0} + \binom{n}{1} = n+1 < \frac{n(n-1)}{2} = \binom{n}{2} < \binom{n}{k}.
\end{equation*}
Therefore we can write
\begin{equation*}
2\sum_{s=0}^{k-1}\binom{n}{s} = 2\left[\binom{n}{0} + \binom{n}{1} + \sum_{s=2}^{k-1}\binom{n}{s}\right] 
\ls 2\left[\binom{n}{k} + \sum_{s=k+1}^{2k-2}\binom{n}{s}\right],
\end{equation*}
and this gives us
\begin{align*}
4\sum_{s=0}^{k-1}\binom{n}{s} &\ls 2\sum_{s=0}^{k-1}\binom{n}{s} + 2\sum_{s=k}^{2k-2}\binom{n}{s} 
\\
&=2^n - 2\binom{n}{2k-1} - \binom{n}{2k} < 2^n - 2\binom{n}{2k},
\end{align*}
which implies the desired inequality.
\item[$\bm{n=4k+2}$ for some $\bm{k\gs 2}$.] Then $\lfloor\frac{n-1}{4} \rfloor = k\gs 2$, and we also have that
\begin{equation*}
2^n = 2 \sum_{s=0}^{2k}\binom{n}{s} \ +\ \binom{n}{2k+1} \ = \ 2\sum_{s=0}^k\binom{n}{s} + 2\sum_{s=k+1}^{2k}\binom{n}{s} \ +\ \binom{n}{2k+1}. 
\end{equation*}
Again combining the facts that $2\binom{n}{2k} > \binom{n}{2k+1}$ and $\binom{n}{0} + \binom{n}{1} + \binom{n}{2} < \binom{n}{3}\ls \binom{n}{k+1}$, we can conclude that
\begin{align*}
4\sum_{s=0}^k\binom{n}{s} &\ls 2\sum_{s=0}^k\binom{n}{s} + 2\sum_{s=k+1}^{2k-1}\binom{n}{s} 
\\
&=2^n - 2\binom{n}{2k} - \binom{n}{2k+1} < 2^n - 2\binom{n}{2k+1},
\end{align*}
as before.
\item[$\bm{n=4k+1}$ for some $\bm{k\gs 2}$.] Then $\lfloor\frac{n-1}{4}\rfloor = k\gs 2$, and we also have that
\begin{equation*}
2^n = 2 \sum_{s=0}^{2k}\binom{n}{s}\ = \ 2\sum_{s=0}^k\binom{n}{s} + 2\sum_{s=k+1}^{2k}\binom{n}{s}. 
\end{equation*}
By the fact that $\binom{n}{0} + \binom{n}{1} + \binom{n}{2} < \binom{n}{3}\ls \binom{n}{k+1}$, we can conclude that
\begin{equation*}
4\sum_{s=0}^k\binom{n}{s}\ls 2\sum_{s=0}^k\binom{n}{s} + 2\sum_{s=k+1}^{2k-1}\binom{n}{s}
= 2^n - 2\binom{n}{2k}.
\end{equation*}
\item[$\bm{n=4k+3}$ for some $\bm{k\gs 2}$.] Then $\lfloor\frac{n-1}{4}\rfloor = k\gs 2$, and we also have that
\begin{equation*}
2^n = 2 \sum_{s=0}^{2k+1}\binom{n}{s}\ = \ 2\sum_{s=0}^k\binom{n}{s} + 2\sum_{s=k+1}^{2k+1}\binom{n}{s}. 
\end{equation*}
By the fact that $\binom{n}{0} + \binom{n}{1} < \binom{n}{2} < \binom{n}{k+1}$, we can conclude that
\begin{equation*}
4\sum_{s=0}^k\binom{n}{s}\ls 2\sum_{s=0}^k\binom{n}{s} + 2\sum_{s=k+1}^{2k}\binom{n}{s}
= 2^n - 2\binom{n}{2k+1}.
\end{equation*}
\end{description}
Finally, regarding the asymptotic behaviour, we can observe the following:
\begin{equation*}
\binom{n}{\lfloor \frac{n-1}{4}\rfloor} \ls \sum_{s=0}^{\lfloor (n-1)/4\rfloor}\binom{n}{s} \ls \frac{n}{2}\cdot \binom{n}{\lfloor \frac{n-1}{4}\rfloor}.
\end{equation*}
Then the desired estimate follows from Stirling's formula.
\end{proof}

\begin{theorem}\label{thm:illum-small-unit-subcubes}
Let $n\gs 4$. We can find a subset ${\cal T}^n$ of $G^n(1/(n+1))$ (the `big' set in Proposition \ref{prop:general-illuminating-set}, when we set $\delta=\frac{1}{n+1}$) which has cardinality
\begin{equation*}
\abs{{\cal T}^n} = 4\sum_{s=0}^{\lfloor (n-2)/4\rfloor}\binom{n-1}{s} \ls 2^{n-1}-4,
\end{equation*}
and the following property: let $\B\in \SU^n$ with $m_\B \ls \frac{n}{2}$; then $\B$ is illuminated by the set
\begin{equation*}
\bigl[\I^{n-1}\bigl(\tfrac{1}{n}\bigr)\times\{0\}\bigr]\cup {\cal T}^n.
\end{equation*}
It follows that, for any such body $\B$, $\,\I(\B) \ls 2^n - 4$ (and in sufficiently high dimensions, $\I(\B) \ls 2^{n-1} + \left(\frac{256}{27}\right)^{\frac{n}{4}\,+\,O(\log n)}$).
\end{theorem}
\begin{proof}
Let $\B\in \SU^n$ with $m_\B\ls \frac{n}{2}$, and let $x\in \partial \B$. We distinguish two cases.
\begin{itemize}
\item If $|x_n| < 1$, then we invoke Corolloary \ref{cor:affine-set} with the affine set $\{\xi\in \R^n: \xi_n = x_n\}$ and Theorem \ref{thm:illum-everything} to deduce that one of the directions in $\bigl[\I^{n-1}\bigl(\tfrac{1}{n}\bigr)\times\{0\}\bigr]$ will illuminate $x$.
\item If $|x_n| =1$, then $\|x\|_\infty = 1$ and therefore, by the main assumption for $\B$, we must have $|M_x|\ls \frac{n}{2}$. Moreover, $n\in M_x$, and thus $|M_x\cap [n-1]|\ls \frac{n}{2}-1$. We know by Lemma \ref{lem:trickier-deep-illumination} that it suffices to find a direction $d_x$ from $G^n(1/(n+1))$ which deep illuminates ${\rm Proj}_{M_x}(x)$ (the projection of $x$ onto the subspace spanned by $e_k, k\in M_x$). We agree from the beginning that the maximum (in absolute value) coordinate of $d_x$ will be the $n$-th one, so for the remaining coordinates we just have to pick the signs carefully so that, if $k\in M_x\cap [n-1]$, then $\sign(x_k)=-\sign(d_{x,k})$. Let $M_x^+=\{k\ls n: |x_k| = \|x\|_\infty\ \hbox{and}\ x_k>0\}$ and $M_x^-=M_x\setminus M_x^+$. It suffices to pick the direction so that 
\begin{equation}\label{eqp1:thm:illum-small-unit-subcubes}
\hbox{$d_{x,k}=-\frac{1}{n+1}$ for all $k\in M_x^+\cap [n-1]$ and $d_{x,k}=+\frac{1}{n+1}$ for all $[n-1]\setminus M_x^+$,}
\end{equation}
or alternatively so that
\begin{equation}\label{eqp2:thm:illum-small-unit-subcubes}
\hbox{$d_{x,k}=+\frac{1}{n+1}$ for all $k\in M_x^-\cap [n-1]$ and $d_{x,k}=-\frac{1}{n+1}$ for all $[n-1]\setminus M_x^-$.}
\end{equation}
We can achieve this by going over all subsets $C_r$ of $[n-1]$ of a certain cardinality $r$, and by setting all entries with index in $C_r$ equal to $-\frac{1}{n+1}$ (say) and all remaining entries equal to $+\frac{1}{n+1}$ (the negatives of the directions we construct in this way will be those directions which have $+\frac{1}{n+1}$ in entries with index in $C_r$, and $-\frac{1}{n+1}$ in all remaining entries; in other words, from a direction which satisfies \eqref{eqp1:thm:illum-small-unit-subcubes}, we can pass to a direction which satisfies \eqref{eqp2:thm:illum-small-unit-subcubes} by taking negatives). 

To be more efficient in how many directions we will include in ${\cal T}^n$, we also note that either $\abs{M_x^+\cap [n-1]}\ls \abs{M_x^-\cap [n-1]}$ or the reverse inequality is true. Hence, we can choose a direction $d_x$ for $x$ which satisfies either \eqref{eqp1:thm:illum-small-unit-subcubes} or \eqref{eqp2:thm:illum-small-unit-subcubes} based on which subset is smaller. In either case, we can check that
\begin{equation*}
\min\{\abs{M_x^+\cap [n-1]}, \abs{M_x^-\cap [n-1]}\} \ls \frac{n-2}{4},
\end{equation*}
and therefore we can set
\begin{align*}
{\cal T}^n:=&\bigl[\tfrac{1}{n+1}\cdot\big\{(\epsilon_1,\epsilon_2,\ldots, \epsilon_{n-1}): \epsilon_i = -1 \ \hbox{for at most $\tfrac{n-2}{4}$ indices $i\in [n-1]$}\bigr\}\bigr]\times\{\pm 1\}
\\
&\ \, \bigcup\  \bigl[\tfrac{1}{n+1}\cdot\big\{(\epsilon_1,\epsilon_2,\ldots, \epsilon_{n-1}): \epsilon_i = +1 \ \hbox{for at most $\tfrac{n-2}{4}$ indices $i\in [n-1]$}\bigr\}\bigr]\times\{\pm 1\}.
\end{align*}  
It follows that
\begin{equation*}
\abs{{\cal T}^n} = 2\cdot 2\sum_{s=0}^{\lfloor (n-2)/4\rfloor}\binom{n-1}{s},
\end{equation*}
and, as explained above, some direction from ${\cal T}^n$ illuminates the point $x$.
\end{itemize}
We have thus shown that $\bigl[\I^{n-1}\bigl(\tfrac{1}{n}\bigr)\times\{0\}\bigr]\cup {\cal T}^n$ illuminates all boundary points of $\B$. It remains to recall Lemma \ref{lem:subset-cardinality}, which gives the desired bounds on the cardinality of ${\cal T}^n$.
\end{proof}

\begin{remark}\bfit{Using pairs of opposite directions.}
In all the results so far, wherever we had to use the set $\I^n(1/(n+1))$ or $\I^{n-1}(1/n)$, we could select it to be any of the sets whose existence Proposition \ref{prop:general-illuminating-set} guarantees.

Of course we can more specifically choose it to be the set that we constructed by the 2nd method in the proof of Proposition \ref{prop:general-illuminating-set}, the `combinatorial construction'. In that case, we will also have that \emph{it is formed from pairs of opposite directions}. This in turn will also imply that the full illuminating sets that we work with in Theorems \ref{thm:illum-no-unit-squares} and \ref{thm:illum-small-unit-subcubes} have this property too.
\end{remark}

In the next subsection we finally come to a setting in which we no longer have the flexibility to pick any eligible construction for $\I^n$: we will specifically have to work with the `combinatorial construction' that we gave.

\subsection{General case: large unit subcubes allowed}

The following theorem deals with the general case, and should be viewed as one of the two main results of this paper: we can now show, for all $n\gs 3$, that, if $\B\in \SU^n$ satisfies ${\rm dist}(\B,[-1,1]^n) > 1$, then $\II(\B) \ls 2^n-2$, and we can do so for small or big $m_\B\in \{1,2,\ldots, n-1\}$ simultaneously and in a unified way. The only thing that may require some adaptations, and which can affect the values of the parameters that we will use, is how close ${\rm dist}(\B, [-1,1]^n)$ is to 1.

\begin{theorem}\label{thm:illum-everything-but-the-cube}
Let $n\gs 3$, $\B\in\SU^n$ and suppose that ${\rm dist}(\B,[-1,1]^n) > 1$ (in other words, $m_\B < n$). 

Let $\I^n\bigl(\frac{1}{n+1}\bigr)$ be the `combinatorially constructed' set from Proposition \ref{prop:general-illuminating-set}, that is, the subset of $G^n\bigl(\frac{1}{n+1}\bigr)$ that we define using the 2nd method in the proof.

\medskip

We can find some sufficiently small $\eta > 0$ (which will depend only on the dimension $n$ and on how small $[{\rm dist}(\B,[-1,1]^n) -1]$ is) such that the set
\begin{gather*}
\Bigl[\I^n\bigl(\tfrac{1}{n+1}\bigr) \backslash\Bigl\{\pm\bigl(+1,+\tfrac{1}{n+1},+\tfrac{1}{n+1},\ldots,+\tfrac{1}{n+1},+\tfrac{1}{n+1},+\tfrac{1}{n+1}\bigr), \  \phantom{\I^n\bigl(\frac{1}{n+1}\bigr)}
\\
\phantom{\pm\bigl(+1,+\tfrac{1}{n+1},\ldots,+\tfrac{1}{n+1},+\tfrac{1}{n+1}\bigr)} \pm\Bigl(+\tfrac{1}{n+1},+\tfrac{1}{n+1},+\tfrac{1}{n+1},\ldots,+\tfrac{1}{n+1},{\bf -\!\!\!\!\!-\frac{1}{n+1}},+1\Bigr)\Bigr\}\Bigr]  \phantom{\I^n\bigl(\frac{1}{n+1}\bigr)}
\\
\bigcup \ \,\Bigl\{\pm\bigl(+1,+\tfrac{1}{n+1},+\tfrac{1}{n+1},\ldots,+\tfrac{1}{n+1},\,{\bm \eta},\, +\tfrac{1}{n+1}\bigr)\Bigr\}
\end{gather*}
will illuminate $\B$.

We can conclude that $\I(\B) \ls 2^n-2$.
\end{theorem}

The key ingredient of the proof to the theorem is the following lemma:
\begin{lemma}\label{lem:deep-illum-with-one-zero}
Let $\delta\in (0,1)$, and let $n\gs 3$. Suppose $\I^n(\delta)$ is the `combinatorially constructed' set from Proposition \ref{prop:general-illuminating-set}, and consider its subset
\begin{equation*}
\I^n_{-2}(\delta):= \I^n(\delta)\setminus\bigl\{\pm(+\delta,+\delta,\ldots,+\delta,+\delta,{\bm{-\delta}},+1)\bigr\}.
\end{equation*}
Then $\I^n_{-2}(\delta)$ can deep illuminate every vector in $\R^n\setminus\{\vec{0}\}$ which has at least one zero coordinate.
\end{lemma}
\begin{proof}
Recall that the set $\I^n(\delta)$ given by the 2nd method in the proof of Proposition \ref{prop:general-illuminating-set} is constructed by dividing it into four parts:
\begin{itemize}
\item the first part consists of all the directions $d$ which satisfy $\sign(d_{n-2})=\sign(d_{n-1})=\sign(d_n)$. Moreover, the maximum (in absolute value) coordinate of $d$ is one of the first $n-2$.
\item The second part consists of all the directions $d$ which satisfy $\sign(d_{n-2})=-\sign(d_{n-1})=-\sign(d_n)$. Moreover, the maximum (in absolute value) coordinate of $d$ is the $(n-1)$-th one.
\item The third part consists of all the directions $d$ which satisfy $\sign(d_{n-2})=\sign(d_{n-1})=-\sign(d_n)$. Moreover, the maximum (in absolute value) coordinate of $d$ is the $n$-th one.
\item Finally, the fourth part of $\I^n(\delta)$ consists of all directions $d$ which satisfy $\sign(d_{n-2})=-\sign(d_{n-1})=\sign(d_n)$. Moreover, the maximum (in absolute value) coordinate of $d$ is the $n$-th one. Observe that this is the part from which we remove two directions.
\end{itemize}
Let $x\in \R^n\setminus\{\vec{0}\}$, and assume $x_{i_0}=0$ for some $i_0\in [n]$. Assume initially that $\{n-2,n-1,n\}\cap \cZ_x\neq \emptyset$ (note that this covers all possibilities when $n=3$). 

If $n\in \cZ_x$ (in which case we can set $i_0=n$), we just need to deep illuminate ${\rm Proj}_{e_n^\perp}(x)$ by a direction $d^\prime\in \I^{n-1}(\delta)$, and then pick the corresponding direction $d= (d^\prime, \sign(d^\prime_{n-1})\delta)\in \I^n(\delta)$ (note that this direction will be contained in the first two parts of $\I^n(\delta)$, hence it is also in $\I^n_{-2}(\delta)$).

If instead $n\notin \cZ_x$, but $n-2\in \cZ_x$ (in which case we could choose $i_0=n-2$), we distinguish two further subcases: 
\begin{itemize}
\item[-] if $n-1\in \cZ_x$ as well, then we can deep illuminate $x$ using the directions from the third part of $\I^n(\delta)$ (since this part contains all possible combinations of signs for the first $n-3$ coordinates and the $n$-th coordinate, and moreover the maximum (in absolute value) coordinate is the $n$-th one, for which we have $x_n\neq 0$).
\item[-] If $n-1\notin \cZ_x$, then we can use directions from either the second part or the third part of $\I^n(\delta)$, based on whether $\sign(x_{n-1})=\sign(x_n)$ or not (we use the second part of $\I^n(\delta)$ in the former case, and the third part in the latter case).
\end{itemize}

Next assume that $\{n-2,\,n\}\cap \cZ_x=\emptyset$. If we still have $n-1\in \cZ_x$ (in which case we can set $i_0=n-1$), then we can use directions from either the first part or the third part of $\I^n(\delta)$ to deep illuminate $x$, based on whether $\sign(x_{n-2})=\sign(x_n)$ or not (we use the first part of $\I^n(\delta)$ in the former case (in fact we first deep illuminate ${\rm Proj}_{[e_{n-1},e_n]^\perp}(x)$ by a direction $d^\prime \in \I^{n-2}(\delta)$, and then work with the direction $d=(d^\prime,\sign(d^\prime_{n-2})\delta,\sign(d^\prime_{n-2})\delta)\in \I^n(\delta)$), while in the latter case we use the third part).

\medskip

Observe at this point that we have completed the proof in the case that $n=3$. For the rest of the proof assume that $n\gs 4$, and note that the only possibility that remains is when $\{n-2,n-1,n\}\cap \cZ_x=\emptyset$ and $\cZ_x\cap [n-3]\neq \emptyset$. 
If it is NOT the case that $\sign(x_{n-2})=-\sign(x_{n-1})=\sign(x_n)$, then we definitely have to use a direction from the first three parts of $\I^n(\delta)$ to deep illuminate $x$; fortunately all these directions remain in $\I^n_{-2}(\delta)$.

\smallskip
 
If instead $\sign(x_{n-2})=-\sign(x_{n-1})=\sign(x_n)$ and $x_{i_0}=0$ for some $i_0\in [n-3]$, then we pick a direction $d$ from the fourth part of $\I^n(\delta)$ that satisfies $\sign(d_i)=-\sign(x_i)$ whenever $x_i\neq 0$, and also satisfies $\sign(d_j)=-\sign(d_n)=-d_n$ for all $j\in \cZ_x\subset [n-3]$ (note that we have at least one such index here, the index $i_0$). This (unique, by our definition) direction deep illuminates $x$, and moreover it is not one of the directions $\pm(+\delta,+\delta,\ldots,+\delta,+\delta,-\delta,+1)$ (since these directions satisfy $\sign(d^\prime_j)=\sign(d^\prime_n)$ for all $j\in [n-3]$). Thus $d$ belongs to $\I^n_{-2}(\delta)$ too.

\smallskip

This completes the proof of the lemma.
\end{proof}

\noindent \emph{Proof of Theorem \ref{thm:illum-everything-but-the-cube}.} Let 
\begin{equation*}
\theta_\B:=\|e_1+e_2+\cdots+e_{n-1}+e_n\|_\B^{-1} = \|{\bm 1}\|_\B^{-1}.
\end{equation*}
We have that $\theta_\B < 1$ since $\B$ is not the cube $[-1,1]^n$. In fact, as already observed in Subsection \ref{subsec:Tikh-and-our-method}, it holds that $\theta_\B=\frac{1}{{\rm dist}(\B,\,[-1,1]^n)}$. Pick 
\begin{equation*}
\eta_\B = \frac{1}{n+2}(1-\theta_\B) < \frac{1}{n+1}
\end{equation*}
and write $\I^n_{-2}\bigl(\frac{1}{n+1},\eta_\B\bigr)$ for the set in the statement of the theorem which one gets if $\eta=\eta_\B$.

\bigskip

Let $x\in \partial \B$, and assume first that $|M_x|\ls n-1$. Then by Lemma \ref{lem:trickier-deep-illumination} (and also Remark \ref{rem:generalise-deep-illumination}) we know that, if a direction $d\in \I_{-2}^n\big(\frac{1}{n+1},\eta_\B\bigr)$ deep illuminates the projection ${\rm Proj}_{M_x}(x)$ of $x$ onto the coordinates with indices in $M_x$, then $d$ will also illuminate $x$. Note that $\I_{-2}^n\big(\frac{1}{n+1},\eta_\B\bigr)$ and $\I_{-2}^n\bigl(\frac{1}{n+1}\bigr)$ only differ in two elements of them, and for those two directions of $\I^n_{-2}\bigl(\frac{1}{n+1}\bigr)$ that have been replaced, we have only changed one of their small coordinates to something even smaller in absolute value (without changing the sign). Thus, by Lemma \ref{lem:deep-illum-with-one-zero}, $\I_{-2}^n\big(\frac{1}{n+1},\eta_\B\bigr)$ will also deep illuminate all the non-zero vectors of $\R^n$ which have at least one zero coordinate, and the projections ${\rm Proj}_{M_x}$ of the boundary points that we consider here do have this property.

\medskip

Consider now the remaining boundary points of $\B$: these are precisely the point $\theta_B {\bm 1}$
and all its coordinate reflections. Let $y$ be one of these points. If we can find a direction $d\in \I_{-2}^n\big(\frac{1}{n+1},\eta_\B\bigr)$ which deep illuminates $y$, then we know that $d$ also illuminates it (in this case we need $\sign(d_i)=-\sign(y_i)$ for all $i\in [n]$ since $y$ has no zero coordinates). In this way, we can take care of all coordinate reflections of $\theta_B {\bm 1}$ except for the following two:
\begin{equation*}
y_{\pm} = \pm (\theta_\B, \theta_\B,\ldots,\theta_\B,{\bm -\theta_\B}, \theta_\B).
\end{equation*}
Consider 
\begin{equation*}
y_+ = (\theta_\B, \theta_\B,\ldots,\theta_\B,{\bm -\theta_\B}, \theta_\B) 
\end{equation*}
and the direction
\begin{equation*}
d_1 = \bigl(-1, -\tfrac{1}{n+1},-\tfrac{1}{n+1},\ldots, -\tfrac{1}{n+1},-\eta_\B,-\tfrac{1}{n+1}\bigr).
\end{equation*}
Then
\begin{equation*}
y_+ + \theta_\B d_1 = \bigl(0,\, (1-\tfrac{1}{n+1})\theta_\B,\ldots, (1-\tfrac{1}{n+1})\theta_\B,\  -(1+\eta_\B)\theta_\B,\  (1-\tfrac{1}{n+1})\theta_\B\bigr).
\end{equation*}
Compare the coordinates of this vector with those of the following convex combination in $\B$:
\begin{align*}
w_0&:=\bigl(1-\tfrac{1}{n+2}\bigr)\theta_\B{\bm 1} + \tfrac{1}{n+2}e_{n-1} 
\\
&\phantom{:}= \bigl((1-\tfrac{1}{n+2})\theta_\B, (1-\tfrac{1}{n+2})\theta_\B,\ldots,(1-\tfrac{1}{n+2})\theta_\B, \ \,\theta_\B + \tfrac{1}{n+2}(1-\theta_\B),\ \, (1-\tfrac{1}{n+2})\theta_\B\bigr).
\end{align*}
For all $i\in [n]\setminus\{n-1\}$, we have that the $i$-th coordinate of $y_++\theta_\B d_1$ is non-negative and strictly smaller than the $i$-th coordinate of $w_0$. Moreover,
\begin{equation*}
|(y_+ + \theta_\B d_1)_{n-1}| = \theta_\B + \eta_\B\theta_\B < \theta_\B + \eta_\B = w_{0,n-1}.
\end{equation*}
Therefore, we can apply Lemma \ref{lem:smaller-coordinates} and conclude that $y_+ + \theta_\B d_1\in \intr\B$.

Similarly we illuminate the boundary point $y_-= -y_+$ using the direction $-d_1$.\qed

\begin{remark}\bfit{Using pairs of opposite directions: Reprise.} Theorem \ref{thm:illum-everything-but-the-cube} shows that, even in the general case, we can use illuminating sets which consist of pairs of opposite directions. Thus, as we also said in the Introduction, in the class of 1-symmetric convex bodies we answer an additional conjecture by Lassak, who surmised in \cite{Lassak-1984} that this should be possible for all origin-symmetric bodies.

This also allows us to settle the X-ray conjecture by Bezdek and Zamfirescu for all 1-symmetric convex bodies. Recall that the X-ray number $X(K)$ of a convex body $K$ in $\R^n$, as proposed by Soltan, is the minimum number $M$ of non-zero vectors $u_1, u_2,\ldots, u_M$ such that, for every $p\in \partial K$, we will have $(p+\R u_i)\cap \intr K \neq \emptyset$ for some $i\in [M]$. Clearly, $X(K)\ls \II(K)\ls 2X(K)$, while Bezdek and Zamfirescu \cite{Bezdek-Zamfirescu-1994} conjectured that we must have $X(K)\ls 3\cdot 2^{n-2}$ for all $K\subset \R^n$ (see \cite{Bezdek-Zamfirescu-1994}, and also e.g. \cite{Bezdek-Kiss-2009}, for further details on this conjecture). 

Of course, in our setting, given that all our illuminating sets can be chosen to consist of pairs of opposite directions, we obtain that $X(\B)\ls 2^{n-1}$ for all $\B\in \SU^n$, $n\gs 3$ (and in fact $X(\B)\ls 2^{n-1} - 1$ if $\B$ is not the cube).

It should also be noted that, in sufficiently high dimensions, the X-ray conjecture for 1-symmetric bodies followed by Tikhomirov's method too. Indeed, recall that Tikhomirov's method gives the following (see Remark \ref{rem:Tikh-proof-gap}): every $B\in \SU^n$ which satisfies ${\rm dist}(\B, [-1,1]^n) > 1$ can be illuminated by one of the sets $T_1$ and $T_3$ (using the notation of Propositions B and C). One can quickly check that the set $T_1$ consists of $2^{n-1}-2$ pairs of opposite directions and 3 additional directions whose negative is not included. Thus if $\B$ is illuminated by $T_1$, then $X(\B) \ls (2^{n-1} -2) + 3$ (by keeping only one direction per pair, along with all 3 directions which are not in a pair). 
On the other hand, if $\B$ is illuminated by $T_3$, and if we also assume that $n$ is sufficiently large, then (by Tikhomirov's bounds on $\abs{T_3}$) we will have that $X(\B)\ls \II(\B) \ls \abs{T_3}\ls 2^{n-1}+\frac{2^n}{n} < 3\cdot 2^{n-2}$.
\end{remark}

\section{Illumination for $(1+\delta)$-symmetric convex bodies}\label{sec:almost-1-sym}

We can use the results in the previous sections to settle the Illumination Conjecture for a slightly wider class of convex bodies, which comes up as a natural enlargement of the class of 1-symmetric convex bodies (we thank Vladimir Troitsky for posing this question to us).

\begin{definition}
Let $K\subset {\mathbb R}^n$ be a convex body, and let $\alpha\gs 1$. We say that $K$ is $\alpha$-symmetric if, for every $x\in K$, for every choice of signs $\epsilon_i\in \{\pm 1\}, i\in [n]$, and for every permutation $\pi\in S_n$, we have that
\begin{equation*}
(\epsilon_1 x_{\pi(1)}, \epsilon_2 x_{\pi(2)},\ldots, \epsilon_{n-1}x_{\pi(n-1)},\epsilon_n x_{\pi(n)})\in \alpha K.
\end{equation*}
\end{definition}

One simple consequence of the definition is that an $\alpha$-symmetric convex body will contain the origin in its interior.
If the parameter $\alpha$ is really close to 1 (how close this needs to be will become clear later), then, as we will see, $K$ will have sufficiently small geometric distance to a genuine 1-symmetric convex body, and will thus satisfy the Illumination Conjecture too.

\begin{lemma}\label{lem:alpha-sym-to-1-sym}
Let $\alpha > 1$, and let $K\subset {\mathbb R}^n$ be an $\alpha$-symmetric convex body. Then we can find a $1$-symmetric convex body $\B_0$ such that 
\begin{equation*}
\B_0 \subset K\subset \alpha \B_0
\end{equation*}
(equivalently, such that ${\rm dist}(K,\B_0)\ls \alpha$).
\end{lemma}
\begin{proof}
Since $K$ contains the origin in its interior, we can consider the following intersection, which will be non-empty:
\begin{equation*}
\B_0 := \bigcap_{\epsilon_i\in \{\pm 1\},\, \pi\in S_n}[{\rm diag}(\epsilon_1,\epsilon_2,\ldots,\epsilon_n)\cdot P_\pi] (K)
\end{equation*}
where $P_\pi$ is the permutation matrix corresponding to the permutation $\pi$, ${\rm diag}(\epsilon_1,\epsilon_2,\ldots,\epsilon_n)$ is the diagonal matrix with diagonal entries $\epsilon_1, \epsilon_2,\ldots,\epsilon_n$, and where we consider all possible ways of combining a choice of signs $\epsilon_i\in \{\pm 1\}, i\in [n]$ and a permutation $\pi\in S_n$ (including the one giving us the identity matrix). 

Then $\B_0$ is a convex body, and moreover it is 1-symmetric. Indeed, let $x\in \B_0$. Then 
for every choice of signs $\eta_i\in \{\pm 1\}, i\in [n]$, and for every permutation $\sigma\in S_n$, we will have
\begin{align*}
(\eta_1 x_{\sigma(1)}, \eta_2 x_{\sigma(2)},\ldots, \eta_n x_{\sigma(n)})&\in 
\\
 [{\rm diag}(\eta_1,\eta_2,\ldots,\eta_n)&\cdot P_\sigma]\left(\bigcap_{\epsilon_i\in \{\pm 1\},\, \pi\in S_n}[{\rm diag}(\epsilon_1,\epsilon_2,\ldots,\epsilon_n)\cdot P_\pi] (K)\right) =
\\
\bigcap_{\epsilon_i\in \{\pm 1\},\, \pi\in S_n}\!\!\!\bigl([{\rm diag}(\eta_1,\eta_2,\ldots,\eta_n)&\cdot P_\sigma]\cdot[{\rm diag}(\epsilon_1,\epsilon_2,\ldots,\epsilon_n)\cdot P_\pi]\bigr) (K)
\\
\intertext{because ${\rm diag}(\eta_1,\eta_2,\ldots,\eta_n)\cdot P_\sigma$ represents an invertible transformation. We can further see that this last intersection is}
&=\bigcap_{\epsilon_i\in \{\pm 1\}, \,\pi\in S_n}\!\!\![{\rm diag}(\eta_1\epsilon_{\sigma(1)},\eta_2\epsilon_{\sigma(2)},\ldots,\eta_n\epsilon_{\sigma(n)})\cdot P_{\sigma\circ\pi}](K)\quad 
\\
&=\  \B_0
\end{align*}
since ${\rm diag}(\eta_1\epsilon_{\sigma(1)},\eta_2\epsilon_{\sigma(2)},\ldots,\eta_n\epsilon_{\sigma(n)})\cdot P_{\sigma\circ\pi}$ runs over all possible combinations of a choice of signs $\widetilde{\epsilon}_i\in \{\pm 1\}, i\in [n]$ and a permutation $\widetilde{\pi}\in S_n$ if $\epsilon_i\in \{\pm 1\}, i\in [n]$ and $\pi\in S_n$ do the same.

\medskip

In addition, $\B_0$ satisfies the required inclusions. The first inclusion is immediate. To confirm the second one, let $x\in K$. Then $x\in \alpha K$ too (since $K$ contains the origin and it is convex, while $\alpha > 1$). Furthermore, for any non-trivial combination of a choice of signs $\epsilon_i\in \{\pm 1\}, i\in [n]$ and of a permutation $\pi\in S_n$, our main assumption for $K$ gives that
\begin{gather*}
(\epsilon_1 x_{\pi(1)}, \epsilon_2 x_{\pi(2)},\ldots, \epsilon_{n-1}x_{\pi(n-1)},\epsilon_n x_{\pi(n)})\in \alpha K
\\[0.5em]
\Rightarrow\ \ x\in \alpha [P_{\pi^{-1}}\cdot {\rm diag}(\epsilon_1,\epsilon_2,\ldots,\epsilon_n)] (K) = 
\alpha [{\rm diag}(\epsilon_{\pi^{-1}(1)},\epsilon_{\pi^{-1}(2)},\ldots,\epsilon_{\pi^{-1}(n)})\cdot P_{\pi^{-1}}](K).
\end{gather*}
Since ${\rm diag}(\epsilon_{\pi^{-1}(1)},\epsilon_{\pi^{-1}(2)},\ldots,\epsilon_{\pi^{-1}(n)})\cdot P_{\pi^{-1}}$ runs over all possible combinations of a choice of signs $\widetilde{\epsilon}_i\in \{\pm 1\}, i\in [n]$ and a permutation $\widetilde{\pi}\in S_n$ if $\epsilon_i\in \{\pm 1\}, i\in [n]$ and $\pi\in S_n$ do the same, we conclude that 
\begin{equation*}
x\in \alpha \left[\bigcap_{\widetilde{\epsilon}_i\in \{\pm 1\}, \,\widetilde{\pi}\in S_n}[{\rm diag}(\widetilde{\epsilon}_1,\widetilde{\epsilon}_2,\ldots,\widetilde{\epsilon}_n)\cdot P_{\widetilde{\pi}}] (K)\right] = \alpha \B_0.
\end{equation*}
The proof is complete.
\end{proof}

Next we recall a semicontinuity argument by Naszódi \cite{Naszodi-2009} and adapt it to our situation. First we need to recall how we pass from illumination to covering.

\begin{remark}\label{rem:illum2cov}
Let $K\subset {\mathbb R}^n$ be a convex body, and let ${\cal D}=\{d_1,d_2,\ldots,d_M\}$ be a set of directions which illuminates $K$ (obviously ${\mathfrak I}(K) \ls M$). Then for every $x\in \partial K$ we can find $d_x\in {\cal D}$ and $\varepsilon_x>0$ such that
\begin{equation*}
x+\varepsilon_xd_x \in \intr(K) \quad \Leftrightarrow\quad x\in - \varepsilon_xd_x + \intr(K).
\end{equation*}
Clearly we can also do that if $y\in \intr(K)$: $y\in - \varepsilon_yd_y + \intr(K)$ for some $d_y\in {\cal D}$ and $\varepsilon_y>0$. Thus
\begin{equation*}
K\subset \bigcup_{z\in K} \bigl(-\varepsilon_zd_z + \intr(K)\bigr)
\end{equation*} 
and by compactness, for each $1\ls i\ls M$ we can find $\varepsilon_1^i, \varepsilon_2^i, \ldots,\varepsilon_{s_i}^i > 0$ such that
\begin{equation*}
K\subset \cup_{i=1}^M \cup_{j=1}^{s_i} \bigl(-\varepsilon_j^id_i + \intr(K)\bigr).
\end{equation*}
It remains to note that, because of convexity, for any point $x\in K$, the containment $x\in -\varepsilon_j^id_i + \intr(K)$ implies also that $x\in -(\min_{1\ls j\ls s_i}\varepsilon_j^i)d_i + \intr(K)$ (indeed, if $\varepsilon_0^i=\min_{1\ls j\ls s_i}\varepsilon_j^i$, and $\varepsilon_{j_x}^i > \varepsilon_0^i$, then $x+\varepsilon_{j_x}^i d_i\in \intr(K)$ implies that $\bigl(1-\frac{\varepsilon_0^i}{\varepsilon_{j_x}^i}\bigr)x + \frac{\varepsilon_0^i}{\varepsilon_{j_x}^i}(x+\varepsilon_{j_x}^id_i)\in \intr(K)$ too).

Thus we have found $\tau_1,\tau_2,\ldots, \tau_M\in {\mathbb R}^n\setminus\{0\}$ such that
\begin{equation*}
K\subset \cup_{i=1}^M (\tau_i +\intr(K)).
\end{equation*}
\end{remark}

\medskip

By compactness again, we can also find $0<\lambda_0<1$ such that $K\subset \cup_{i=1}^M (\tau_i +\intr(\lambda_0 K))\subset \cup_{i=1}^M (\tau_i +\lambda_0 K)$ (we can first assume for convenience that $K$ has been translated to contain the origin in its interior).

We can now observe (see e.g. \cite[Proposition 2.2]{Naszodi-2009}) that, if $L$ is another convex body which satisfies ${\rm dist}(K,L)\ls 1 + \frac{1-\lambda_0}{2}$, then we can also write
\begin{equation*}
L\subset \cup_{i=1}^M (\widetilde{\tau}_i + \widetilde{\lambda_0}\, L)
\end{equation*}
for some $\widetilde{\lambda_0} \ls \lambda_0 + \lambda_0\frac{1-\lambda_0}{2} < 1$. Thus $\II(L)\ls M$ too. In fact, because the illumination number is an affine invariant, we can even allow translates when considering distance, and more generally use the Banach-Mazur distance: as long as we can find $a,b\in \R^n$ and $T\in {\rm GL}(n)$ such that
\begin{equation*}
T(L)-b \subseteq K-a\subseteq \bigl(1 + \tfrac{1-\lambda_0}{2}\bigr)(T(L)-b),
\end{equation*}
the above conclusion holds (compare also with Lemma \ref{lem:small-distance-to-cube} from Section \ref{sec:prelims}).

Naszódi in \cite{Naszodi-2009} uses this to conclude that, if $\sup_{K\subset {\mathbb R^n}}\II(K) \ls A$, there is a common homothety factor $\lambda_n = \lambda_{n,A} \in (0,1)$ with the property that, for any convex body $K\subset {\mathbb R}^n$, $N(K,\lambda_nK)\ls A$ (recall that $N(Q_1,Q_2)$ stands for the covering number of a set $Q_1$ by another set $Q_2$).
The key ingredient here is the compactness of the Banach-Mazur compactum, so we can have the same conclusion within smaller compact classes of convex bodies in ${\mathbb R}^n$.

\begin{proposition}\label{prop:cover-by-homothets}
For every $n\gs 3$, we can find $\lambda_{1,n}\in (0,1)$ such that, for every convex body $\B\subset {\mathbb R}^n$ affinely equivalent to a $1$-symmetric convex body, we will have $N(\B, \lambda_{1,n}\B) \ls 2^n$.

\smallskip

Furthermore, for every $\delta_0\in (0,1)$, we also have the following: we can find $\lambda_{2,n} = \lambda_{2,n,\delta_0}\in (0,1)$ such that, if $\B$ is affinely equivalent to a $1$-symmetric convex body $\widetilde{\B}$ which satisfies ${\rm dist}(\widetilde{\B},\,[-1,1]^n) \gs 1+\delta_0$, then $N(\B, \lambda_{2,n}\B)\ls 2^n-2$.
\end{proposition}
\begin{proof}
We essentially reproduce Naszódi's argument, but adapted to the situation at hand (and with some details added for clarity). The illumination number, which coincides with $N(K,\intr K)$, is an affine invariant, and in fact so are the quantities $N(K,\mu K)$ for any $\mu > 0$. Thus, we can assume that every convex body we consider is in $\SU^n$. By Theorem \ref{thm:illum-everything}, we know that, for every $\B\in \SU^n$, $\II(\B)\ls 2^n$. Then, as explained before, by compactness we can find some $\mu_\B\in (0,1)$ such that $N(\B, \mu_\B \B) \ls 2^n$ (we can choose $\mu_\B$ in a deterministic way too, by e.g. setting it equal to $\frac{1+\mu_{\B, {\rm inf}}}{2}$ where $\mu_{\B, {\rm inf}}$ is the infimum of all homothety factors $\mu\in (0,1)$ such that $N(\B,\mu \B)\ls 2^n$).

Assume towards a contradiction that $\sup_{\B\in \SU^n} \mu_\B = 1$ (note that this is equivalent to assuming that $\sup_{\B\in \SU^n}\mu_{\B, {\rm inf}} = 1$), and consider a sequence of convex bodies $\B_i\in \SU^n$ such that $\mu_{\B_i}\to 1$. Recall that, by our normalisation, for all these bodies we have
\begin{equation*}
CP_1^n = {\rm conv}\{\pm e_i: i\in [n]\}\subseteq \B_i \subseteq [-1,1]^n,
\end{equation*}
and thus we can apply Blaschke's selection theorem to obtain a subsequence $(\B_{i_j})$ of $(\B_i)_{i\in {\mathbb N}}$ and a convex body $K_0$ that $(\B_{i_j})$ converges to in the Hausdorff metric. Since every $\B_{i_j}$ is 1-symmetric, and in fact in $\SU^n$, we can conclude that $K_0$ is also in $\SU^n$. Thus we have $\II(K_0)\ls 2^n$ and $N(K_0,\mu_{K_0} K_0)\ls 2^n$ for some $\mu_{K_0} < 1$.

By what we said above, if $L$ is any other convex body whose Banach-Mazur distance ${\rm dist}_{BM}(K_0,L)$ from $K_0$ is $\ls 1 + \frac{1-\mu_{K_0}}{2}$, then we will be able to write $N(L,\mu_1 L) \ls 2^n$ with $\mu_1$ some constant $< \frac{\mu_{K_0}+1}{2} < 1$.

In particular, this will be true for the convex bodies $\B_{i_j}$ once $j$ gets sufficiently large (since small Hausdorff distance from $K_0$ will imply ${\rm dist}_{BM}(K_0,\B_{i_j})$ gets sufficiently small too). But this gives $\limsup\mu_{\B_{i_j},{\rm inf}}\ls \mu_1 < 1$, contradicting our assumption about the sequence $(\B_i)_{i\in {\mathbb N}}$. 

Therefore, we must have $\sup_{\B\in \SU^n} \mu_\B < 1$, and thus we can set $\lambda_{1,n}$ equal to this supremum.

\bigskip

For the second part of the statement, we use a completely analogous argument, except that now we restrict our attention to convex bodies $\widetilde{\B}$ in $\SU^n$ such that 
\begin{equation*}
\|e_1+e_2+\cdots+e_n\|_{\widetilde{\B}}^{-1} = \bigl({\rm dist}(\widetilde{\B},\,[-1,1]^n)\bigr)^{-1} \ls \frac{1}{1+\delta_0}.
\end{equation*}
For such convex bodies, we have shown in Theorem \ref{thm:illum-everything-but-the-cube} that ${\mathfrak I}(\widetilde{\B})\ls 2^n-2$. Moreover, for any sequence $(\widetilde{\B}_i)_{i\in {\mathbb N}}$ of such bodies, Blaschke's selection theorem will give us a convergent subsequence which converges to a convex body $\widetilde{K}_0$ which also satisfies
\begin{equation*}
\|e_1+e_2+\cdots+e_n\|_{\widetilde{K}_0}^{-1}  \ls \frac{1}{1+\delta_0}.
\end{equation*}
Hence, a straightforward adjustment of the previous argument gives the desired conclusion. 
\end{proof}

We are now ready to verify the Illumination Conjecture for $(1+\delta_n)$-symmetric convex bodies too, where $\delta_n$ will be a small positive constant (allowed to depend on the dimension). In fact, the main idea in the following argument can also be used to settle the X-ray conjecture for such convex bodies (see the remark at the end)

\begin{theorem}
For every $n\gs 3$ there is $\alpha_n > 1$ such that, if $K\subset {\mathbb R}^n$ is an $\alpha_n$-symmetric convex body (or is affinely equivalent to such a body), then $K$ satisfies the Illumination Conjecture. In other words, $\II(K)\ls 2^n-1$ unless $K$ is a parallelepiped, in which case $\II(K) = 2^n$.
\end{theorem}
\begin{proof}
As mentioned in the Introduction, Livshyts and Tikhomirov \cite{Livshyts-Tikhomirov-2020} have shown that, for every $n\gs 3$, there exists a small constant $\widetilde{\delta}_n$ such that, if $K$ is not a parallelepiped and satisfies ${\rm dist}_{BM}(K,[-1,1]^n) \ls 1+\widetilde{\delta}_n$, then $\II(K)\ls 2^n-1$ (in fact, the upper bound here is optimal: as they show in their paper, there are convex bodies, as close to the $n$-dimensional cube as one wants, which have illumination number equal to $2^n-1$). We will eventually choose $\alpha_n \ls \sqrt{1+\widetilde{\delta}_n}$. 

At the same time, based on the previous proposition, we can find $\lambda_2 = \lambda_{2,n,\widetilde{\delta}_n}\in (0,1)$ such that, for all $\B\in \SU^n$ with ${\rm dist}(\B,[-1,1]^n)\gs \sqrt{1+\widetilde{\delta}_n}$, we will have $N(\B, \lambda_2 \B) \ls 2^n-2$. But then, for any convex body $L$ which satisfies ${\rm dist}_{BM}(L, \B)\ls 1 + \frac{1-\lambda_2}{2} $ for some $\B$ as above, we will be able to deduce that $\II(L)\ls 2^n-2$.

\medskip

Fix now
\begin{equation*}
\alpha_n = \min\left\{\sqrt{1+\widetilde{\delta}_n},\ 1 + \frac{1-\lambda_2}{2}\right\}
\end{equation*}
and consider an $\alpha_n$-symmetric convex body $K$ in ${\mathbb R}^n$ which is not a parallelepiped. By Lemma \ref{lem:alpha-sym-to-1-sym} we can find a 1-symmetric convex body $\B_0$ such that ${\rm dist}(K,\B_0)\ls \alpha_n$. We consider two possibilities.
\begin{itemize}
\item ${\rm dist}(\B_0,[-1,1]^n) < \sqrt{1+\widetilde{\delta}_n}$. Then 
\begin{equation*}
{\rm dist}(K,[-1,1]^n)\ls {\rm dist}(K,\B_0)\cdot {\rm dist}(\B_0,[-1,1]^n) < \alpha_n\cdot \sqrt{1+\widetilde{\delta}_n}\ls 1+\widetilde{\delta}_n,
\end{equation*}
and thus by \cite{Livshyts-Tikhomirov-2020} we obtain that $\II(K)\ls 2^n-1$.
\item ${\rm dist}(\B_0,[-1,1]^n) \gs \sqrt{1+\widetilde{\delta}_n}$. Then, given how we chose $\lambda_2$ above, we will have that
\begin{equation*}
{\rm dist}(K,\B_0)\ls \alpha_n \ls 1 + \frac{1-\lambda_2}{2}
\end{equation*}
implies that $\II(K)\ls 2^n-2$.
\end{itemize}
The proof is complete.
\end{proof}

\begin{remark}
Going back to Remark \ref{rem:illum2cov} and Proposition \ref{prop:cover-by-homothets}, let us examine the arguments there more carefully: our purpose for doing this is to confirm that we have also settled (through the above) the X-ray conjecture for $\alpha_n$-symmetric convex bodies, if $\alpha_n >1$ is small enough (depending on the dimension). 

Indeed, we can check that, if ${\cal D}= \{d_1,d_2,\ldots,d_M\}$ is an illuminating set for a convex body $K\subset \R^n$ (where $M\gs \II(K)$), then we can write 
\begin{equation*}
K\subset \cup_{i=1}^M (\tau_i +\intr(K)),
\end{equation*}
with the translation vectors $\tau_i$ satisfying:
\begin{equation*}
\bigl\{\tau_i: 1\ls i\ls M\bigr\} = \bigl\{-\widehat{\varepsilon}_0d_i: 1\ls i\ls M\bigr\}
\end{equation*}
for some small enough $\widehat{\varepsilon}_0>0$. Afterwards, by compactness, we can also find $\lambda_0\in (0,1)$ so that $K\subset \cup_{i=1}^M (-\widehat{\varepsilon}_0d_i +\lambda_0 K)$.
But then, if $L_1$ is another convex body in $\R^n$ satisfying the precise inclusions
\begin{equation*}
L_1 \subseteq K \subseteq \bigl(1+\tfrac{1-\lambda_0}{2}\bigr)L_1,
\end{equation*}
we will be able to write
\begin{equation*}
L_1 \subset K \subset \cup_{i=1}^M (-\widehat{\varepsilon}_0d_i +\lambda_0 K) \subset \cup_{i=1}^M \left(-\widehat{\varepsilon}_0d_i +\lambda_0 \bigl(1+\tfrac{1-\lambda_0}{2}\bigr)L_1\right),
\end{equation*}
which shows that, for every boundary point $y$ of $L_1$, 
\begin{equation*}
\hbox{there is $i=i_y\in [M]$ so that:}\quad
y\in -\widehat{\varepsilon}_0d_i + \intr L_1 \quad \Leftrightarrow\quad 
y+\widehat{\varepsilon}_0d_i \in \intr L_1.
\end{equation*}
In other words, ${\cal D}$ is an illuminating set for $L_1$ as well.

Now, going more carefully over the proof of Proposition \ref{prop:cover-by-homothets}, we see that for every $n\gs 3$ there exists $\lambda_{1,n}\in (0,1)$ so that, if $\B\in \SU^n$, then we will be able to find some $\varepsilon_\B > 0$ such that
\begin{equation*}
\B \subseteq \!\!\!\!\bigcup_{d\in \I^n(1/(n+1))} \!\!\!\!\bigl(-\varepsilon_\B d + \lambda_{1,n}\B\bigr),
\end{equation*}
where $\I^n(1/(n+1))$ is any of the sets satisfying Proposition \ref{prop:general-illuminating-set} (that we have fixed beforehand). In other words, now we also keep track of what the translation vectors will be (up to the scaling factor $\varepsilon_\B$) when we cover $\B$ by smaller homothetic copies of itself. 

We finally also assume that {\bf we have chosen a construction of $\I^n(1/(n+1))$ which consists of pairs of opposite directions} (e.g. the construction that we have already presented).

But then, if we set $\widetilde{\alpha}_n = 1 + \frac{1-\lambda_{1,n}}{2}$, and we consider $K\subset \R^n$ which is affinely equivalent to an $\widetilde{\alpha}_n$-symmetric convex body, we can conclude that some affine image of $K$ is illuminated by the set $\I^n(1/(n+1))$, hence by a set which consists of $2^{n-1}$ pairs of opposite directions. Thus, both this affine image and $K$ itself satisfy the X-ray conjecture (and also give an affirmative answer to Lassak's question about `efficient' illumination using pairs of opposite directions; this is true even though $K$ may not necessarily have a centre of symmetry).
\end{remark}


\noindent {\bf Acknowledgements.} The first-named author would like to thank Matt Yan for his invaluable assistance in writing a program that enabled him to develop intuition for the discovery of some of the examples/counterexamples described in Section \ref{sec:counterexamples}. 

Part of writing up the final version of this paper was done while the two authors were in residence at the Hausdorff Research Institute for Mathematics for the programme ``Synergies between modern probability, geometric analysis and stochastic geometry''. The authors are grateful to the institute and the organisers for the hospitality and the excellent working conditions. The second-named author is partially supported by an NSERC Discovery Grant.

\section*{Appendix A}

Here we justify Claim F from the end of Section \ref{sec:counterexamples} (which recovers \cite[Proposition 5]{Tikhomirov-2017}).

\smallskip

\emph{Proof of Claim F.} Fix $n\gs 3$ and $\B\in \SU^n$ with $1\neq {\rm dist}(\B,[-1,1]^n) < 2$. Given the normalisation of $\B$, and since $1 < {\rm dist}(\B,[-1,1]^n) < 2$, we can find $\epsilon_0\in (0,1/2)$ such that
\begin{equation}\label{eqp:ClaimF-1}
\bigl(\tfrac{1}{2}+\epsilon_0\bigr){\bm 1} =\bigl(\tfrac{1}{2}+\epsilon_0\bigr) \sum_{i=1}^n e_i\in \B
\end{equation}
(we might as well assume here that $\epsilon_0$ is maximum possible).

Let us suppose that $\B$ is not illuminated by the set $T_1$. Then, as the proof (and not just the statement) of \cite[Lemma 7]{Tikhomirov-2017} shows, \underline{$\B$ must contain} a boundary point of the form 
\begin{equation}\label{eqp:ClaimF-2}
x_0:=(\gamma_\B,\gamma_\B,\ldots,\gamma_\B,\gamma_\B, 1)
\end{equation}
with $\gamma_\B \in (0,1)$ and such that 
\begin{equation}\label{eqp:ClaimF-3}
\B\subset \{y\in \R^n: y_i+y_n \ls x_{0,i} + x_{0,n} = \gamma_\B + 1\}
\end{equation}
for any $i\in [n-1]$; in other words the vectors $e_i + e_n$, $i\in [n-1]$, are outer normals of $\B$ at $x_0$. 

\medskip

\emph{We justify this claim first:} due to \cite[Lemma 6]{Tikhomirov-2017}, and given the existence of the point in \eqref{eqp:ClaimF-1} (which makes the inequality in the statement of \cite[Lemma 6]{Tikhomirov-2017} impossible to satisfy), the only points in $\partial\B$ that wouldn't be illuminated right away by some direction in $T_1$ are points $x$ which have only non-zero entries and such that $x_j < 0$ for all $j\ls n-1$.

\smallskip

Moreover, if $|x_n| < 1$, then we could also use Corollary \ref{cor:affine-set} to confirm that the direction $e_1+e_2+\cdots+e_{n-1}$ illuminates a point $x$ as above. Thus, we must also assume that $|x_n|=1$.

\medskip

Finally, we can rule out the possibility that such a point $x$, which is not illuminated by any direction in $T_1$, could satisfy $\min_{j\in [n-1]}|x_j| < \max_{j\in [n-1]}|x_j|$: indeed, if we had $|x_{j_1}|=\min_{j\in [n-1]}|x_j|$ and $|x_{j_2}| > |x_{j_1}|$ for some $j_1,j_2\in [n-1]$, then $x$ would be illuminated by the direction
\begin{equation*}
-e_{j_1}-\sign(x_n)e_n+\sum_{i\in [n-1]\setminus\{j_1\}} e_i
\end{equation*}
in $T_1$. This is because such a boundary point $x$ would not have any outer normal vectors $v$ with $v_{j_1}\neq 0$ and $v_j=0$ for all $j\in [n-1]\setminus\{j_1\}$ (because any such vector $v$ would satisfy $\langle P_{j_1,j_2}(x), v\rangle > \langle x,v\rangle$, where $P_{j_1,j_2}$ is the transformation swapping the $j_1$-th and the $j_2$-th coordinates). Recall now \cite[Lemma 4]{Tikhomirov-2017} which tells us that, if $w$ is an arbitrary outer normal at $x$ with $w_{j_1}\neq 0$, then we will have that $|w_{j^\prime}|\gs |w_{j_1}| > 0$ for both $j^\prime=j_2$ and $j^\prime=n$; moreover, recall the basic observation that $w_i\cdot x_i\gs 0$ for all $i\in [n]$ (following from the 1-unconditionality). Thus 
\begin{equation*}
\left\langle -e_{j_1}-\sign(x_n)e_n+\sum_{i\in [n-1]\setminus\{j_1\}} e_i,\  w\right\rangle = |w_{j_1}| - |w_n| + \sum_{i\in [n-1]\setminus\{j_1\}}(-|w_i|) \ <\  0
\end{equation*}
since $|w_{j_1}|-|w_n|\ls 0$ and since there will be at least one more non-zero summand, the term $-|w_{j_2}|$, in the expansion of the above dot product. Clearly, it is even simpler to see that the corresponding dot product will be negative if $w_{j_1}=0$. Thus the desired conclusion follows from Criterion A.

\bigskip

From the above discussion, we are guaranteed that the `problematic' boundary points of $\B$ have the form
\begin{equation*}
\tilde{x}\equiv \tilde{x}_{\gamma,\pm} =(-\gamma,-\gamma,\ldots,-\gamma,-\gamma, \pm 1)
\end{equation*}
for some $\gamma > 0$. We can also assume that $\gamma=\gamma_\B$ is maximum possible (given that, for smaller values of $\gamma$, the corresponding points will not be extreme, so they can be handled quickly once we illuminate all extreme (or possibly extreme) points of $\B$). 

Recall now that ${\rm dist}(\B, [-1,1]^n) > 1$, and thus $\gamma_\B < 1$ and $|\tilde{x}_j| < |\tilde{x}_n|$ for all $j\ls n-1$. But then we can go over the above discussion again, as appropriate, to obtain that a direction of the form
\begin{equation*}
-e_{j_1}-\sign(x_n)e_n+\sum_{i\in [n-1]\setminus\{j_1\}} e_i
\end{equation*}
for some $j_1\in [n-1]$ would illuminate $\tilde{x}\,$ \underline{\bf except in the case that} $\,\tilde{x}$ would have an outer normal $v_0$ of the form $|v_{0,j_1}|=|v_{0,n}| \neq 0$, and $v_{0,i}=0$ for all $i\in [n]\setminus\{j_1,n\}$. In the latter case we would get the condition
\begin{equation*}
\gamma_\B + 1 = |\tilde{x}_{j_1}| + |\tilde{x}_n| = \max\{|y_{j_1}| + |y_n| : y\in \B\},
\end{equation*}
which is what we claim at the beginning.

\bigskip

For the remainder of this proof, we rely on the assumptions in \eqref{eqp:ClaimF-1} and in the paragraph containing \eqref{eqp:ClaimF-2}-\eqref{eqp:ClaimF-3}. Given that $1+\gamma_\B < 2$, and that $1+\gamma_\B = \max\{|y_i|+|y_n|: y\in \B\}$, where $i$ is any index $\ls n-1$, we cannot have $e_i+e_n\in \B$. In other words, the type of 1-symmetric convex body $\B$ that we have now focused on satisfies $\|e_i+e_j\|_\B > 1$ for all $i\neq j$ in addition to the other conditions; let us set $\beta_0 = \|e_i+e_j\|_\B^{-1}$.

\smallskip

It remains to observe that, because the point in \eqref{eqp:ClaimF-1} has been chosen to be a boundary point of $\B$, we cannot have all the coordinates of $x_0$ be bigger than $\frac{1}{2}+\epsilon_0$. This shows that
\begin{equation}\label{eqp:ClaimF-4}
\gamma_\B \ls \frac{1}{2}+\epsilon_0.
\end{equation}

\medskip

We can now check that $T_2$ illuminates $\B$ by arguing as follows: if $x\in \partial\B$ satisfies $|x_n|=1$, then $|x_n|$ will be the unique maximum (in absolute value) coordinate of $x$, and hence $\pm e_n$ will illuminate $x$ according to Lemma \ref{lem:unique-max-coord} (the same conclusion will hold if $|x_n| > \beta_0$, as this will still imply that $|x_n|$ is the unique maximum coordinate of $x$).

\medskip

On the other hand, if we have $|x_n| < 1$, we could invoke Corollary \ref{cor:affine-set}, and try to illuminate ${\rm Proj}_{e_n^\perp}(x)$ as a boundary point of $P_{x_n}\equiv{\rm Proj}_{e_n^\perp}(\B_{x_n})$ where $\B_{x_n}=\{y\in \B: y_n= x_n\}$: we would like to show that $P_{x_n}$ (viewed as a subset of $\R^{n-1}$) is illuminated by $\{-1,1\}^{n-1}$. Recalling Lemma \ref{lem:small-distance-to-cube}, it would suffice to prove that $P_{x_n}$ contains 
\begin{equation*}
\tfrac{\|e_1\|_{P_{x_n}}^{-1}}{2}(e_1+e_2+\cdots+e_{n-1})
\end{equation*} 
in its interior.

Note that it is easy to conclude this if $|x_n| \ls \frac{1}{2}+\epsilon_0$, given that $\|e_1\|_{P_{x_n}}^{-1} \ls 1$ in all cases, and moreover here we will have 
\begin{equation*}
\bigl(\tfrac{1}{2}+\epsilon_0, \tfrac{1}{2}+\epsilon_0,\ldots,\tfrac{1}{2}+\epsilon_0,\tfrac{1}{2}+\epsilon_0,\,x_n\bigr)\in\B_{x_n}.
\end{equation*}

\smallskip

Thus it remains to consider the regime $\frac{1}{2}+\epsilon_0 < |x_n| < 1$. Let us try to estimate $\|e_1\|_{P_{x_n}}$ in such cases. Consider a point $y$ on the relative boundary of $\B_{x_n}$ that has only one additional non-zero coordinate (besides the $n$-th one); 
let's say $y_i$ is the non-zero coordinate.
Then we will have $|y_i|+|y_n| = |y_i| + |x_n| \ls 1+\gamma_\B$, which shows that $\|e_1\|_{P_{x_n}}^{-1} = |y_i| \ls 1+\gamma_\B - |x_n|$.

Next we look for $\lambda_{x_n}\in (0,1)$ so that
\begin{equation*}
\lambda_{x_n} + (1-\lambda_{x_n})\bigl(\tfrac{1}{2}+\epsilon_0\bigr) = |x_n| \quad \Leftrightarrow \quad \lambda_{x_n}=\frac{|x_n| - \frac{1}{2}-\epsilon_0}{\frac{1}{2}-\epsilon_0}\,.
\end{equation*}
We then have that the convex combination
\begin{equation*}
\lambda_{x_n}(\gamma_\B,\gamma_\B,\ldots,\gamma_\B,\gamma_\B, \sign(x_n)) + (1-\lambda_{x_n})\bigl(\tfrac{1}{2}+\epsilon_0, \tfrac{1}{2}+\epsilon_0,\ldots,\tfrac{1}{2}+\epsilon_0,\tfrac{1}{2}+\epsilon_0,\sign(x_n)\bigl(\tfrac{1}{2}+\epsilon_0\bigr)\bigr)
\end{equation*}
will be contained in $\B_{x_n}$. It holds that
\begin{align*}
\lambda_{x_n}\gamma_\B + (1-\lambda_{x_n})\bigl(\tfrac{1}{2}+\epsilon_0\bigr) &= \frac{|x_n| - \frac{1}{2}-\epsilon_0}{\frac{1}{2}-\epsilon_0}\gamma_\B \ +\  \frac{1-|x_n|}{\frac{1}{2}-\epsilon_0}\bigl(\tfrac{1}{2}+\epsilon_0\bigr)
\\[0.3em]
&= \frac{|x_n|\gamma_\B\,+\,\bigl(\frac{1}{2}+\epsilon_0\bigr)(1-\gamma_\B-|x_n|)}{\frac{1}{2}-\epsilon_0}\,.
\end{align*}
We would like this expression to be
\begin{equation*}
> \frac{1+\gamma_\B-|x_n|}{2}.
\end{equation*}
But this is equivalent to
\begin{equation*}
 \frac{1-|x_n|}{2} + 3\epsilon_0(1-|x_n|) + 2|x_n|\gamma_\B \  >\  \frac{3}{2}\gamma_\B +\epsilon_0\gamma_\B,
\end{equation*}
which we can rewrite as
\begin{equation*}
\frac{1-|x_n|}{2} + 3\epsilon_0(1-|x_n|) + 2\bigl(|x_n|-\tfrac{1}{2}-\epsilon_0\bigr)\gamma_\B \  >\  \frac{\gamma_\B}{2} -\epsilon_0\gamma_\B,
\end{equation*}
and furthermore as
\begin{align*}
\bigl(\tfrac{1}{2}-\epsilon_0\bigr)\bigl(\lambda_{x_n}2\gamma_\B + (1-\lambda_{x_n})&\tfrac{1}{2}\bigr) + 3\epsilon_0(1-|x_n|) \ >\ \frac{\gamma_\B}{2} -\epsilon_0\gamma_\B
\\
&\Leftrightarrow\quad \bigl(\tfrac{1}{2}-\epsilon_0\bigr)\bigl(\lambda_{x_n}2\gamma_\B + (1-\lambda_{x_n})\tfrac{1}{2} - \gamma_\B\bigr) + 3\epsilon_0(1-|x_n|) \ >\ 0
\\[0.5em]
&\Leftrightarrow\quad \bigl(\tfrac{1}{2}-\epsilon_0\bigr)\bigl(\lambda_{x_n}\gamma_\B + (1-\lambda_{x_n})\bigl(\tfrac{1}{2}-\gamma_\B\bigr)\bigr) + 3\epsilon_0(1-|x_n|) \ >\ 0.
\end{align*}
Because of \eqref{eqp:ClaimF-4}, we have $\tfrac{1}{2}-\gamma_\B \gs -\epsilon_0$ and therefore
\begin{align*}
&\bigl(\tfrac{1}{2}-\epsilon_0\bigr)\bigl(\lambda_{x_n}\gamma_\B + (1-\lambda_{x_n})\bigl(\tfrac{1}{2}-\gamma_\B\bigr)\bigr) + 3\epsilon_0(1-|x_n|)
\\
&\ \ \gs \bigl(\tfrac{1}{2}-\epsilon_0\bigr)\bigl(\lambda_{x_n}\gamma_\B - (1-\lambda_{x_n})\epsilon_0\bigr)
+ 3\epsilon_0(1-|x_n|) 
\\
&\ \ = \bigl(|x_n| - \tfrac{1}{2}-\epsilon_0\bigr)\gamma_\B + 2\epsilon_0(1-|x_n|) \ >\ 0.
\end{align*}
Thus, in all cases where 
$|x_n| < 1$, the section $\B_{x_n}$ contains a point of the form
\begin{equation*}
(\theta_{x_n},\theta_{x_n},\ldots, \theta_{x_n},\theta_{x_n},\,x_n)
\end{equation*}
with $\theta_{x_n} > \frac{1}{2}\max\{\rho\in (0,1]: \rho e_1 + x_n e_n\in \B\}$. This allows us to apply Lemma \ref{lem:small-distance-to-cube} (and Corollary \ref{cor:affine-set}) and illuminate each such section of $\B$ using the directions $\{-1,1\}^{n-1}\times\{0\}$ (more accurately, illuminate in this way each boundary point of $\B$ that belongs to such a section).

This completes the proof of Claim F. \qed

\bigskip
\bigskip
\bigskip

\noindent \textsc{Department of Mathematical and Statistical Sciences,
University of Alberta, CAB 632, Edmonton, AB, Canada T6G 2G1}

\smallskip

\noindent 
{\it E-mail addresses:} \texttt{wrsun@ualberta.ca, vritsiou@ualberta.ca}

\end{document}